\documentclass[12pt]{amsart}

\usepackage{fullpage}
\usepackage{lipsum}
\usepackage{amsfonts}
\usepackage{amsmath,amssymb,amsthm}
\usepackage{graphicx}
\usepackage{wrapfig}
\usepackage{subcaption}
\usepackage{comment}
\usepackage{cleveref}
\usepackage{rlepsf}
\usepackage{placeins}
\usepackage{tikz}
\usepackage{tikz-cd}

\newtheorem{theorem}{Theorem}[section]
\newtheorem{lemma}[theorem]{Lemma}

\newtheorem{question}[theorem]{Question}

\theoremstyle{definition}
\newtheorem{definition}[theorem]{Definition}
\newtheorem{ex}[theorem]{Example}

\theoremstyle{remark}
\newtheorem{remark}[theorem]{Remark}

\def\I{\mathbb I}

\def\R{\mathbb R}

\def\Z{\mathbb Z}

\def\A{\mathcal A}

\def\V{\mathcal V}
\def\E{\mathcal E}

\def\X{\mathcal X}

\def\U{\mathcal U}
\def\oG{\overline G}
\def\ogamma{\overline \gamma}

\newcommand{\diam}{\mathrm{diam}}
\newcommand{\rank}{\mathrm{rank}}

\newcommand{\inj}{\mathrm{inj}}

\newcommand{\wh}{\widehat} 
\newcommand{\wt}{\widetilde}

\newcommand{\PD}{\mathrm{PD}}
\newcommand{\vr}[2]{\mathrm{VR}(#1;#2)}

\title{Geometric Approaches on Persistent Homology}

\author{Henry Adams}
\email{henry.adams@colostate.edu}
\author{Baris Coskunuzer}
\email{baris.coskunuzer@utdallas.edu}

\thanks{The first author is partially supported by the National Science Foundation under grant No.\ 1934725.
The second author is partially supported by a Simons Collaboration Grant.}

\thanks{\emph{Keywords:}
Persistent homology, geometric topology, width, thick-thin decomposition, power filtration, Vietoris--Rips complexes.}

\thanks{\emph{2020 Mathematics Subject Classification:} 55N31, 55U10, 57R19, 62R40.}

\begin{document}

\maketitle

\begin{abstract}
We introduce several geometric notions, including the width of a homology class, to the theory of persistent homology.
These ideas provide geometric interpretations of persistence diagrams.
Indeed, we give quantitative and geometric descriptions of the ``life span'' or ``persistence'' of a homology class.
As a case study, we analyze the power filtration on unweighted graphs, and provide explicit bounds for the life spans of homology classes in persistence diagrams in all dimensions.
\end{abstract}

\section{Introduction}

In this paper, we investigate persistent homology notions by using geometric topology techniques.
In the past decade, topological data analysis (TDA) has grown substantially, and has proven to be quite useful to understand many phenomena described via different types of data \cite{Carlsson2009}.
During this period, the theoretical foundations of TDA have mostly benefited from algebraic topology, as it was born in part in that domain.
Here, we bring new perspectives, ideas, and terminology to help interpret these powerful methods by using its sister field, geometric topology.

In some applications of topology, persistent homology is considered mostly as a black box by data scientists, where it produces features to be applied to the problem at hand.
These features appear as a \emph{persistence diagram} $\PD_k$, summarizing the $k$-dimensional holes that appear.
While persistence diagrams are known to describe these $k$-dimensional holes, and their life spans are interpreted as the ``size'' of these $k$-dimensional holes, a rigorous mathematical definition for the size of these homology classes has not yet been given in all settings.

In this paper, we give explicit interpretations of the outcomes of persistence diagrams in terms of the geometry of the data.
In order to establish these geometric notions, we chose the simplest setup, where the distances change discretely, i.e.\ unweighted graphs with the power filtration.
We give an explicit geometric description of the persistence diagrams in this case.
In other words, we show how the persistence diagram of the power filtration measures the sizes of holes.
While we mainly focus on unweighted graphs with the power filtration as a case study, the techniques we introduce here are general, and they can be adapted to different settings related to persistent homology.

Our main results are as follows:
For any dimension $k\ge 2$, and for any $k$-cycle $\sigma$ with birth $b$ and death $d$ in the $k$-dimensional persistence diagram $(b,d)\in \PD_k$, we ask if there are upper bounds for $d/b$
in terms of the volume $\|\sigma\|$ (\cref{ques:vol2}), and we prove upper bounds for $d/b$ and in terms of the width $\omega(\sigma)$ (\cref{thm:main-width}).
On the other hand, for $k=1$, we give a complete, explicit description of the persistence diagrams in terms of the lengths of a lexicographically shortest basis of $1$-cycles $\{\gamma_i\}$, i.e.,  $\PD_1=\{(\, 1,\left\lceil \tfrac{|\gamma_i|}{3}\right\rceil\,\}$ (\cref{thm-main1}).
Note that in the related settings of (continuous) metric graphs and geodesic spaces, analogous versions of this result about 1-dimensional persistence were proven in \cite{gasparovic2018complete,virk20201}.

While proving these theorems, we introduce several new geometric notions into the setting of persistent homology, e.g.\ \emph{min-max technique, sweepouts, width of a homology class, thick-thin decomposition, injectivity radius, bracelets.}
The width can be interpreted as the size of a homology class, and we establish its relationship with the life spans in persistence diagrams.
In particular, our results show that for a $k$-cycle $\sigma$ with birth and death $(b,d)\in \PD_k$, the ratio $d/b$ corresponds to the size of the $k$-dimensional cavity $\sigma$ as measured by the ``width'' of the homology class $\sigma$.
Our width approach further relates the size of a homology class with Gromov's filling radius of a $k$-manifold \cite{gromov1983filling}; see \cite{lim2020vietoris} and \cref{ssec:interpret} for further discussion.

From the point of view of TDA on graphs, our results show that higher persistence diagrams indeed contain useful information about the graph.
In many applications, researchers only consider the 0- and 1-dimensional persistence diagrams $\PD_0$ and $\PD_1$ because they are computationally cheaper.
Our results and examples indicate that the diagrams $\PD_k$ capture valuable and geometrically interpretable information about the graph's properties also for dimension $k\ge 2$.
See \cref{ssec:tdaongraphs} for further discussion.

The organization of the paper is as follows.
In \cref{sec:background}, we overview the related work in the subject, and  describe the setting for power filtrations of unweighted graphs.
We give a complete description of the 1-dimensional persistence diagram $\PD_1$ in terms of lengths of $1$-cycles in \cref{sec:PD1}.
In \cref{sec:PD2}, we give some interesting examples to motivate our study of higher-dimensional homology.
In \cref{sec:area}, we ask how the life spans of $2$-cycle and $k$-cycles are related to their area and volume.
Here, we introduce the notions of thick-thin decomposition, injectivity radius, bracelet, and the volume of a homology class.
In \cref{sec:width}, we upper bound the life spans of $2$-cycles by a notion called the width.
In this section, we introduce the geometric notions of the min-max technique, sweepouts, and the width of a homology class.
In \cref{sec:higher}, we generalize our width result to all homological dimensions.
Finally, in \cref{sec:remarks}, we give some concluding remarks.

\section{Background}
\label{sec:background}

\subsection{Related Work on Vietoris--Rips Complexes}

If $G$ is a connected graph, then its vertex set can be equipped with the structure of a metric space, where the distance between any two vertices is the (integer) length of the shortest path between them.
The power filtration of a graph that we study in this paper, namely the clique complexes of the graph powers of $G$, is nothing other than the Vietoris--Rips simplicial complex filtration of the vertex set of $G$.
As Vietoris--Rips complexes transform a metric space into a simplicial complex, they were invented by %Leopold
Vietoris to provide a cohomology theory for metric spaces~\cite{Hausmann1995,lefschetz1942algebraic,Vietoris27}.

If $G$ is a Cayley graph of a group (constructed with respect to a chosen set of generators), then the clique complexes of the graph powers of $G$ are the Vietoris--Rips complexes of the group when equipped with the word metric.
Indeed, Vietoris--Rips complexes were used in geometric group theory by %Eliyahu
Rips as a natural way to thicken a space, and to show that torsion-free hyperbolic groups have Eilenberg--MacLane spaces with finitely many cells~\cite{bridson2011metric}.
% We hope our results will be interesting to geometric group theorests as well.

In applied and computational topology, Vietoris--Rips complexes are used to thicken a data set $X$ in order to approximate its underlying shape.
The shape of a dataset is often reflective of important patterns within~\cite{Carlsson2009}.
Indeed, connected components correspond to segments of the data that could be analyzed separately, circles correspond to periodic or recurrent phenomena, and higher-dimensional features represent further structure in the data.
The Vietoris--Rips complex $\vr{X}{r}$ on metric space $X$ contains $X$ as its vertex set, and a finite subset of $X$ as a simplex if its diameter is at most $r$.
Since we do not know a priori how to choose the thickening scale $r$, the idea of persistent homology is to compute the homology of the Vietoris--Rips complex of data set $X$ over a large range of scale parameters $r$ and to trust those topological features which persist.
The \emph{persistence} or \emph{life span} of a feature is typically defined as its death scale minus its birth scale, and sometimes defined as its death scale divided by its birth scale (as we will use here).
It is understood that the persistence of a $k$-dimensional feature is related to its ``geometric size"; in this paper we make this intuitive notion more precise.

The motivation for using Vietoris--Rips complexes in applied topology is a remarkable theorem due to Latschev \cite{Latschev2001}: for $M$ a closed Riemannian manifold, for scale $r$ sufficiently small depending on the curvature of $M$, and for data set $X$ close to $M$ in the Gromov--Hausdorff distance, we have a homotopy equivalence $\vr{X}{r} \simeq M$.
This result is an analogue of the Nerve Theorem~\cite{Borsuk1948} for Vietoris--Rips complexes, and it has been expanded upon by the manifold reconstruction results in~\cite{ChazalDeSilvaOudot2014,EdelsbrunnerHarer,niyogi2008finding,Chambers2010}, which also rely on the scale being chosen to be sufficiently small.
But as the main idea of persistence is to allow $r$ to vary, the assumption that scale $r$ is kept sufficiently small typically fails in practice.

Indeed, the situation that data scientists are confronted with is that they are given a data set $X$ noisily sampled from an unknown shape $M$.
Without knowing $M$, they do not how to pick the scale parameter $r$ small enough for the above reconstruction guarantees to hold.
As a result, they instead let the scale parameter $r$ in the Vietoris--Rips complexes $\vr{X}{r}$ vary from zero to large.
Hence data scientists construct Vietoris--Rips complexes at large scale parameters on top of their data (and there is efficient software designed to do this~\cite{bauer2021ripser}), even though we do not yet have a mathematical understanding of how these simplicial complexes behave at large scales.
The circle is essentially the only manifold $M$ for which the homotopy types of $\vr{M}{r}$ are known at all scale parameters $r$~\cite{AA-VRS1}, and its proof is built upon approximating the circle via denser and denser graphs~\cite{adamaszek2013clique,AAFPP-J}.
%Persistent homology allows them the ability to vary the scale from small to large, and to encode this multi-resolution view of the data in a concise manner.
Could the homotopy types of Vietoris--Rips complexes of $n$-spheres or other manifolds also be proven via graph approximations?

Applied and computational topology is recently being connected more tightly to quantitative topology, especially the filling radius.
The \emph{filling radius} of a manifold $M$ was used by %Mikhael
Gromov to prove the systolic inequality, which provides a lower bound for the volume of an essential manifold $M$ in terms of the length of the shortest non-contractible loop~\cite{gromov1983filling,gromov2007metric}.
In subsequent work, Katz determined the filling radius of spheres and projective spaces~\cite{katz1983filling,katz1989diameter,katz1991neighborhoods}.
The recent work~\cite{lim2020vietoris,okutan2019persistence} by Lim, M\'{e}moli, and Okutan shows that Vietoris--Rips complexes are strongly connected to quantitative topology: if $M$ is a manifold, then the top-dimensional bar in the persistent homology for the Vietoris--Rips complex filtration has a death time determined by the filling radius.
This same paper proves that any persistent homology bar in the Vietoris--Rips filtration of a metric space $X$ has persistence (birth minus death) upper bounded by the \emph{spread} of $X$.
We instead give bounds on the life span which depend on the individual homology class under consideration.

We view our work, bounding the life spans of the $k$-dimensional homology holes in Vietoris--Rips complexes of unweighted graphs, as being a first step towards injecting ideas from geometric topology into this conversation.
Several key ideas from geometric topology, such as the width and thick-thin decompositions, may not be so well-known to applied topologists.
Nevertheless, in the setting of unweighted graphs, these ideas provide more precise geometric interpretations for the lengths of persistent homology features.
To generalize our techniques to most data analysis settings, we will need to allow non-integer distances: ideas in coarse geometry allowing one to approximate a metric spaces via graphs~\cite{burago2001course, burago2015uniform} may enable generalizations along these lines.

\subsection{Power Filtration for Unweighted Graphs} 
\label{ssec:powerfiltration}

To study the size of a homology class from the perspective of geometric topology, we consider the setting of the power filtration on unweighted graphs, as it provides a discrete and simple setup.
First, we describe some preliminary notions related to homology and to the power filtration on unweighted graphs.

\subsubsection*{Homology} 
For $L$ a simplicial complex and for $k\ge 0$, let $H_k(L)$ denote the simplicial homology of $L$, taken with coefficients in $\Z/2\Z$.
With $\Z/2\Z$ coefficients, we can represent any $k$-chain as a set of simplices, i.e.\ as a sum of simplices where all (nonzero) coefficients are equal to one, which simplifies arguments.
Recall that if $C_k(L)$ is the set of all $k$-chains, then we have a boundary map $\partial \colon C_k(L) \to C_{k-1}(L)$ that satisifes $\partial\circ\partial=0$.
A $k$-cycle is a $k$-chain $\sigma$ satisfying $\partial \sigma = 0$ (which we also may write as $\partial \sigma = \emptyset$ since we are using $\Z/2\Z$ coefficients).

While we use $\Z/2\Z$-coefficients throughout the paper, we expect that many of the ideas can be adapted to other homology coefficients with the necessary modifications made.

\subsubsection*{Metric on $G$} 
Throughout the paper, we assume that $G$ is a finite simple graph that is connected.
For such a graph $G$, let $\V=\{v_i\}$ be the set of vertices in $G$, and let $\E=\{e_{ij}\}$ be the edges in $G$, where $e_{ij}$ represents the edge between the vertices $v_i$ and $v_j$ if it exists.
We define the metric $\rho\colon \V \times \V \to \R$ on the vertex set $\V$ by assigning length $1$ to all edges in $G$.
In particular, if $m_{ij}$ is the smallest number of edges required to get from $v_i$ to $v_j$ in $G$, then the distance between $v_i$ and $v_j$ is defined as $\rho(v_i,v_j)=m_{ij}$.
Notice that as we assumed $G$ is connected, we have $\rho(v_i,v_j)<\infty$ for any vertices $v_i$ and $v_j$.
Let $\diam(G)=\max\{\rho(v_i,v_j)~|~v_i,v_j\in\V\}$ be the diameter of $G$, which we will sometimes denote as $D$.
From the topological data analysis perspective, we use the set of vertices $\V$ as our point cloud, and we use the edges $\E$ to define our metric on this  point cloud.

\subsubsection*{Power Filtration} 
\begin{figure}[h]
	\begin{center}
		$\begin{array}{c@{\hspace{.4in}}c@{\hspace{.4in}}c}
		\relabelbox  {\epsfxsize=1.35in \epsfbox{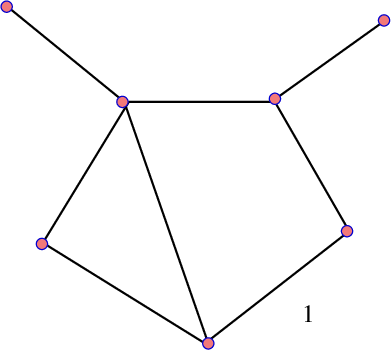}}
		\relabel{1}{ $G$} 				
		\endrelabelbox &
		\relabelbox  {\epsfxsize=1.35in \epsfbox{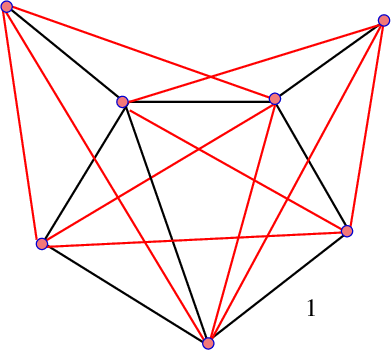}}
		\relabel{1}{ $G^2$} 				
		\endrelabelbox &
		\relabelbox  {\epsfxsize=1.35in \epsfbox{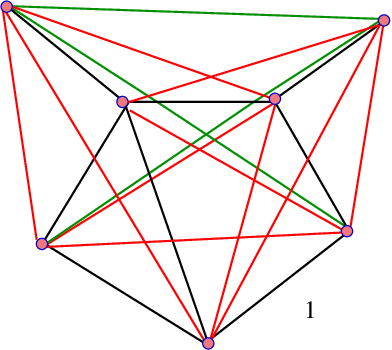}}
		\relabel{1}{ $G^3$} 				
		\endrelabelbox \\	
	\end{array}$
\end{center}
\caption{\footnotesize {\bf Graph Powers.} A graph $G=G^1$ and its graph powers.
	Red edges are added in $G^2$, and green ones are added in $G^3$.
	Note $G^3$ is the complete graph on $7$ vertices since $D=\diam(G)=3$.
	Hence, all higher powers are same, i.e.\ $G^n=G^3$ for $n\ge 3$.}
\label{fig:graph-power2}
\end{figure}

Let $G^n$ be the graph induced by $G$ by adding edges $\{e_{ij}\}$ between vertices $v_i, v_j$ with $\rho(v_i,v_j)\le n$.
In other words, we do not change the nodes of $G$, but if there are vertices at distance $\le n$ in $G$, then in the graph $G^n$ we have an edge between these vertices.
By convention, we set $G^0$ to be the graph with vertex set $\V$ and with no edges.
We will call this new graph $G^n$ the {\em $n^{th}$ power of $G$}; see \cref{fig:graph-power2}.

\setlength{\belowcaptionskip}{-10pt}
	
Let $\wh{G}^n$ be the clique complex of $G^n$.
In more detail, $\wh{G}^n$ contains an $m$-simplex spanning the vertices $v_{i_0}, v_{i_1},\ldots,v_{i_m}$ if for any $i_j,i_k$ with $0\le j,k\le m$, we have $\rho(v_{i_j},v_{i_k})\le n$.
In particular, if there exists a triangle of three edges in $G^n$, then in $\wh{G}^n$ we fill the triangle with a $2$-simplex.
Similarly, if there are $4$ vertices each pairwise connected to each other by edges in $G^n$, then in $\wh{G}^n$ we fill it with a tetrahedron, and so on for simplices of all dimensions.
We remark that $\wh{G}^n$ is the Vietoris--Rips complex of the metric space $\V$ (with the shortest-path metric described above) at scale $n$.
Furthermore, we can safely restrict attention to integer scale parameters, since any shortest path distance in $G$ is a nonnegative integer.
By allowing $n$ to vary, we obtain the following \emph{power filtration} induced by the graph $G$.
\[ \wh{G}^0\subset \wh{G}^1\subset \wh{G}^2\subset \ldots \subset \wh{G}^{D-1} \subset \wh{G}^D\]

Notice that $\wh{G}^0$ is equal to the set of vertices $\V$ in $G$.
Simplicial complex $\wh{G}^1=\wh{G}$ is the clique complex of the original graph $G$.
To form $\wh{G}^2$, we add new edges and cliques to $\wh{G}^1$ accordingly.
In particular, if $\rho(v_i,v_j)=2$, then a new edge $\wt{e}_{ij}$ is added to $\wh{G}^1$.
Similarly, if there is a set of vertices $\{v_{i_0},\ldots,v_{i_k}\}$ where the pairwise distances are at most $n$, then there exists a $k$-simplex $\sigma=[v_{i_0},\ldots,v_{i_k}]$ in $\wh{G}^n$.
Note also that $\wh{G}^n$ is the $(|\V|-1)$-simplex for any $n\ge D=\diam(G)$, and hence is contractible.

For $k\ge 0$, we take $k$-dimensional homology  $H_k$ with coefficients in $\Z/2\Z$.
The points in a persistence diagram represent the birth and death times of a homology class.
In particular, let $\PD_k(G)$ denote the persistence diagram for the $k$-dimensional homology of the power filtration of the graph $G$~\cite{EdelsbrunnerHarer,edelsbrunner2000topological,zomorodian2005computing}.
Then, any persistence diagram point $(b,d)\in \PD_k(G)$ represents a $k$-cycle $\sigma$ that is born in $\wh{G}^b$, and that first becomes homologous to earlier features in $\wh{G}^d$.
In other words, $n=b$ is the birth time for $\sigma$, while $n=d$ is the death time for $\sigma$.

More explicitly, after applying homology with coefficients in the field $\Z/2\Z$ to the power filtration, we obtain the following \emph{persistence module}, i.e.\ the following sequence of vector spaces equipped with linear maps in-between:
\[ H_k(\wh{G}^0)\to H_k(\wh{G}^1)\to H_k(\wh{G}^2)\to \ldots \to H_k(\wh{G}^{D-1}) \to H_k(\wh{G}^D).\]
It follows from~\cite{zomorodian2005computing} or~\cite{ZigzagPersistence,gabriel1972unzerlegbare} that this persistence module decomposes uniquely (up to reordering) as a direct sum 
$\oplus_{i=1}^N \I_{b_i,d_i}$ with $b_i\in\{0,\ldots,D\}$, $d_i\in\{0,\ldots,D\}\cup\{\infty\}$, and $b_i<d_i$.
In this direct sum, each term $\I_{b_i,d_i}$ is an \emph{interval sequence} of the form
\[0\to\ldots\to 0\to\Z/2\Z\to\ldots\to\Z/2\Z\to 0\to\ldots\to 0,\]
where the first copy of the field $\Z/2\Z$ appears in index $b_i$, where the last copy of the field $\Z/2\Z$ appears in index $d_i-1$,\footnote{
In the case when the last copy of the field $\Z/2\Z$ appears in index $D$, then by convention we set $d_i=\infty$.
For power filtrations, since $\wh{G}^D$ is contractible, we obtain only a single persistence diagram point with death value $\infty$, which appears as a single bar of the form $(b,d)=(0,\infty)$ in $\PD_0(G)$.} where all linear maps between adjacent copies of $\Z/2\Z$ are the identity map, and where all other linear maps are the zero map.
Then, the persistence diagram $\PD_k(G)$ is defined as the multiset $\PD_k(G)=\{(b_i,d_i) \mid 1\le i\le N\}$, meaning that we have $N$ different points in the persistence diagram, each of the form $(b_i,d_i)$ for $1\le i\le N$.

The persistence diagram for $0$-dimensional homology with the power filtration, namely $\PD_0(G)$, is easy to understand.
In more detail, while $\rank(H_0(\wh{G}^0))=|\V|$ is the number of vertices, we have $\rank(H_0(\wh{G}^n))=1$ for any $ n\ge 1$ since $G$ is connected.
This means that $\PD_0(G)$ consists of $|\V|-1$ birth-death pairs of the form $(b,d)=(0,1)$, and one birth-death pair of the form $(b,d)=(0,\infty)$.
So $\PD_0(G)=\{(0,1)^{|\V|-1}, (0,\infty)\}$.
Here, $(b,d)^m$ means that the multiplicity of the point $(b,d)$ is equal to $m$, or in other words, that the persistence diagram $\PD_0(G)$ consists of $m$ copies of the point $(b,d)$.

\section{$\PD_1(G)$: Persistence Diagrams in Dimension $1$}
\label{sec:PD1}

Let $G$ be a finite connected graph.
In this section, we give an explicit description of the 1-dimensional persistence diagram $\PD_1(G)$ of the graph $G$ in terms of the lengths of certain loops.
The persistence diagrams $\PD_k(G)$ in homological dimensions $k=0$ and $1$ are the most frequently used in applications, and the most efficient diagrams to compute.

A \emph{path of length $l$} in $G$ is a sequence of vertices $v_0, v_1, \ldots, v_{l-1}, v_l$ such that each $[v_i, v_{i+1}]$ is an edge in $G$ for $0\le i\le l-1$.
If $v_l=v_0$, then this is furthermore a \emph{loop of length $l$} in $G$.

\begin{lemma}
\label{lem:birth-1}
Let $G$ be a finite connected graph.
Every point in $\PD_1(G)$ has birth time $b=1$, and so $\PD_1(G)=\{(1, d_i)\}_i$ for some collection of death times $d_i$.
\end{lemma}

\begin{proof}
This lemma follows from~\cite[Fact~2.1]{adamaszek2013clique}, but we give a stand-alone proof.
Let $(b,d)\in \PD_1(G)$.
Clearly $b\ge 1$, since $\wh{G}^0$ is a disjoint collection of vertices.

We claim that a loop $\alpha$ in $\wh{G}^n$ for $n\ge 1$ is homotopy equivalent in $\wh{G}^n$ to a loop in $G\subseteq \wh{G}^1$.
Let $e=[w^-,w^+]$ be an edge in $\alpha$ with $\rho(w^-,w^+)=t$ for $1\le t\le n$.
If $t=1$, then we leave $e$ unchanged.
Otherwise, let $\{w^-,v_1, v_2,\ldots,v_{t-1},w^+\}$ be the vertices along the shortest path from $w^-$ to $w^+$ in $G$.
Then, $\sigma=[w^-,v_1,\ldots,v_{t-1},w^+]$ is a $t$-simplex in $\wh{G}^n$.
Furthermore, the edge $e=[w^-,w^+]$ is homotopy equivalent in $\sigma$ to the path $\tau=[w^-,v_1]\cup[v_1,v_2]\cup\ldots\cup[v_{t-2},v_{t-1}]\cup[v_{t-1},w^+]$, while fixing the endpoints.
Note that the path $\tau$ is in $\wh{G}$.
Therefore, by applying this process to each edge in $\alpha$, we see that $\alpha$ is homotopy equivalent in $\wh{G}^n$ to a loop in $G\subseteq\wh{G}^1$.
This shows that if $(b,d)$ is a point in $\PD_1(G)$, then we have $b=1$.

Since every point in $\PD_1(G)$ is born at $b=1$, we know that $\PD_1(G)=\{(1, d_i)\}_i$ for some collection of death times $d_i$.
\end{proof}

Let $\lceil.\rceil$ be the ceiling function, i.e.\ $\lceil x \rceil$ is the smallest integer greater than or equal to $x$.

\begin{lemma}
\label{lem:death-3}
A loop $\gamma$ of length $l$ in the finite connected graph $G$ is null-homotopic in $\wh{G}^n$ for $n\ge \lceil\frac{l}{3}\rceil$.
\end{lemma}

\begin{figure}[h]
\begin{center}
		$\begin{array}{c@{\hspace{.4in}}c}

		\relabelbox  {\epsfxsize=2.2in \epsfbox{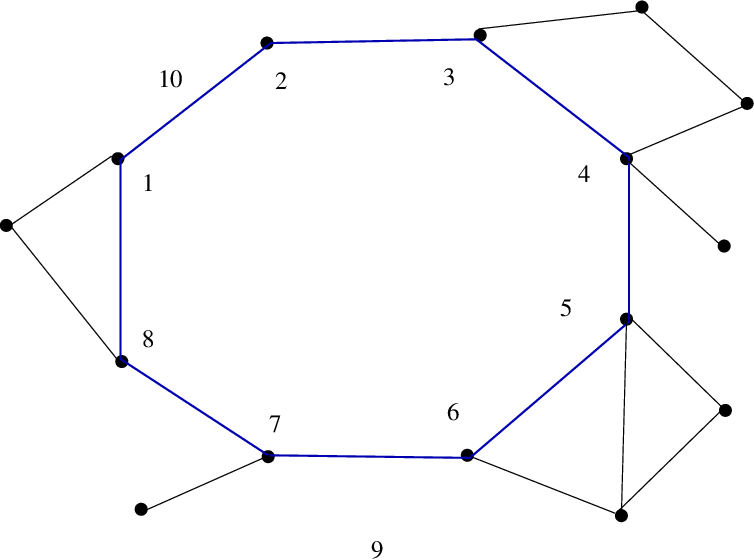}}
		\relabel{1}{\footnotesize $v_0$} 
		\relabel{2}{\footnotesize $v_1$}  
		\relabel{3}{\footnotesize $v_2$} 
		\relabel{4}{\footnotesize $v_3$} 
		\relabel{5}{\footnotesize $v_4$} 
		\relabel{6}{\footnotesize $v_5$} 
		\relabel{7}{\footnotesize $v_6$}  
		\relabel{8}{\footnotesize $v_7$} 
    	\relabel{9}{\small $G$} 
		\relabel{10}{\small $\gamma$} 
		\endrelabelbox &
			
		\relabelbox  {\epsfxsize=2.2in \epsfbox{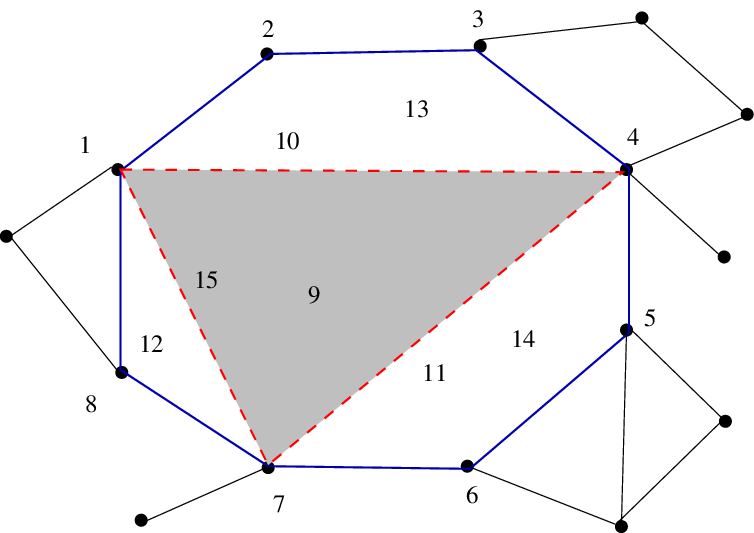}}
		\relabel{1}{\footnotesize $v_0$} 
		\relabel{2}{\footnotesize $v_1$}  
		\relabel{3}{\footnotesize $v_2$} 
		\relabel{4}{\footnotesize $v_3$} 
		\relabel{5}{\footnotesize $v_4$} 
		\relabel{6}{\footnotesize $v_5$} 
		\relabel{7}{\footnotesize $v_6$}  
		\relabel{8}{\footnotesize $v_7$} 
		\relabel{9}{} 
		\relabel{10}{} 
		\relabel{11}{} 
		\relabel{12}{\footnotesize $\Delta_6$} 
		\relabel{13}{\footnotesize $\Delta_0$}
		\relabel{14}{\footnotesize $\Delta_3$}
		\relabel{15}{} 
		\endrelabelbox \\			
	\end{array}$		
\end{center}
\caption{ \footnotesize
\textbf{Spanning Disks.}
On the left, the graph $G$ has a (blue) loop of length $8$.
On the right, this loop is filled by the simplices $\Delta_0$, $\Delta_3$, $\Delta_6$, and the gray triangle.
}
\label{fig_disk}
\end{figure}

\begin{proof}
It suffices to prove the case $n=\lceil\frac{l}{3}\rceil$.
See \cref{fig_disk}.
As $\gamma$ is a loop in $G$ of length $l$, we can write $\gamma=[v_0,v_1]\cup \ldots \cup[v_{l-2},v_{l-1}]\cup [v_{l-1},v_0]$.
Note that $\Delta_0:=[v_0,v_1,\ldots v_n]$, $\Delta_n:=[v_n,v_{n+1},\ldots v_{2n}]$, and $\Delta_{2n}:=[v_{2n},v_{2n+1},\ldots,v_{l-1},v_0]$ are each simplices in the clique complex $\wh{G}^n$.
Furthermore, note that the edge $[v_0,v_n]$ is homotopy equivalent in $\Delta_0$ to the (shortest) arc in $\gamma$ with the same endpoints, and similarly for the edge $[v_n,v_{2n}]$ in $\Delta_n$, and for the edge $[v_{2n},v_0]$ in $\Delta_{2n}$.
Therefore $\gamma$ is homotopy equivalent in $\wh{G}^n$ to the loop $[v_0,v_n]\cup[v_n,v_{2n}]\cup[v_{2n},v_0]$,
which is null-homotopic in the clique complex $\wh{G}^n$ as it is filled by the triangle $[v_0,v_n,v_{2n}]$.
Therefore $\gamma$ is null-homotopic in $\wh{G}^n$.
\end{proof}

Denote the first Betti number of the graph $G$ by $m=\rank(H_1(G))$.
Let the loops $\{\gamma_1,\gamma_2,\ldots,\gamma_m\}$ generate a basis for $H_1(G)$, and let $l_i$ be the length of the loop $\gamma_i$.
We call $\gamma_1,\ldots,\gamma_m$ a {\em lexicographically shortest basis} of $H_1(G)$ if the non-decreasing sequence of lengths $l_1\le l_2\le l_3\le \ldots\le l_m$ is lexicographically smallest among all bases for $H_1(G)$.
See~\cite{gasparovic2018complete,virk20201} for further discussion on this definition.
We are now ready to give the full description of $\PD_1(G)$.

\begin{theorem}
\label{thm-main1}
Let $G$ be a finite connected graph.
Let $\gamma_1,\ldots,\gamma_m$ be a lexicographically shortest basis of $H_1(G)$, the first homology of $G$, and let each $\gamma_i$ have length $l_i$.
Then
\[\PD_1(G)=\{(\, 1,\left\lceil \tfrac{l_i}{3}\right\rceil\, )\mid 1\le i\le m\}.\]
\end{theorem}

Since $G$ is a simple graph, all loops have length at least $3$.
In the lexicographically shortest basis $\gamma_1,\ldots,\gamma_m$, some of the lengths $l_1\le l_2\le \ldots\le l_m$ may be equal to $3$.
If so, then they contribute a persistence diagram point $(b,d)=(1,\left\lceil\tfrac{3}{3}\right\rceil)=(1,1)$ along the diagonal; such points along the diagonal are typically ignored as the death time is equal to the birth time.
Indeed, these loops of length $3$ appear in the power filtration at stage $\wh{G}^1$, when they are immediately filled in since the clique complex contains a 2-simplex filling-in each triangle of three edges.

This theorem should be thought of as an analogue of Theorem~8.10 of~\cite{virk20201} and Theorem~1.1 of~\cite{gasparovic2018complete}, which give a similar characterization for the 1-dimensional persistent homology of Vietoris--Rips or \v{C}ech complexes of metric graphs and geodesics spaces, and indeed we use these results in our proof.
Metric graphs and geodesic spaces are non-discrete metric spaces.
By contrast, our result characterizes the 1-dimensional persistent homology of Vietoris--Rips complexes of the discrete vertex subset of an unweighted graph.

Our proof of \cref{thm-main1} relies not only on \cref{lem:birth-1,lem:death-3}, but also some technical machinery (Vietoris--Rips complexes of geodesic spaces, persistence modules that are indexed over the real numbers instead of over a discrete set, and morphisms between persistence modules~\cite{bauer2014induced,gasparovic2018complete,virk20201}) that is not needed in the rest of our paper.
Therefore, we defer the proof of \cref{thm-main1} to \cref{app:proof}.

\begin{remark}
\cref{thm-main1} can be interpreted as saying that $\PD_1(G)$ detects the number of the essential loops in the clique complex $\wh{G}$, along with their lengths.
In other words, any element $(1,d)\in \PD_1(G)$ with $d>1$ represents that there exists an essential loop of length $\approx 3d$ in the clique complex $\wh{G}$.
\end{remark}

\section{$\PD_k(G)$: Persistence Diagrams in Dimension $2$ and Higher}
\label{sec:PD2}

So far, we have given an explicit description of the $0$- and $1$-dimensional persistent homology for power filtrations, $\PD_0(G)$ and $\PD_1(G)$, in terms of geometric properties of the graph $G$.
In this part, we discuss higher-dimensional persistence diagrams, before providing one-sided generalizations of the $k=0$ and $k=1$ results to higher dimensions in the following sections.

Notice that in our setting, the birth times for all topological features in dimensions $0$ and $1$ are known by construction.
Indeed, all of the $0$-dimensional birth times are $0$ while all the $1$-dimensional birth times are $1$, i.e.\ $PD_0(G)=\{(0,d_j)\}_j$ and $PD_1(G)=\{(1,d_i)\}_i$. 
We start by noting that a direct generalization of \cref{thm-main1} to 2-dimensional persistent homology $\PD_2(G)$ is not true, since not all 2-dimensional features have the same birth scale.
We will show this by using the following interesting examples.

\begin{figure}[h]
\centering
\begin{subfigure}{1in}
\relabelbox  {\epsfxsize=0.9in \epsfbox{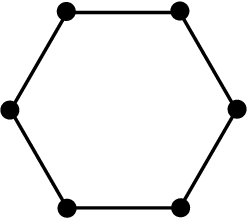}}

\endrelabelbox
\caption{The graph $C_6$}
\end{subfigure}%
\hspace{.2in}
\begin{subfigure}{3in}
\relabelbox  {\epsfxsize=3in \epsfbox{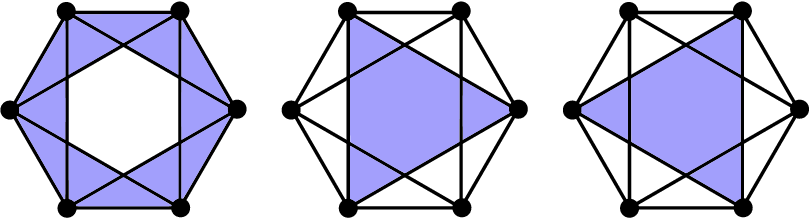}}

\endrelabelbox
\caption{The eight maximal 2-simplices in $\wh{C}_6^2$}
\end{subfigure}
\caption{Power filtration of the graph $C_6$.}
\label{fig:C6}
\end{figure}

\begin{ex}
Let $C_6$ be the cycle of length $6$, namely the graph with 6 vertices and 6 edges, arranged in a loop (\cref{fig:C6}-left).
Since this graph has no triangles, we have $\wh{C}_6^1=C_6$.
Furthermore, the simplicial complex $\wh{C}_6^2$ is homeomorphic to the 2-sphere $S^2$.
Indeed, there are 8 maximal simplices in $\wh{C}_6^2$ (\cref{fig:C6}-right), and they are all 2-simplices.
Six of these 2-simplices, near the boundary of the cycle, glue together to form a cylinder.
The remaining two ``equilateral" 2-simplices in \cref{fig:C6}-right get attached as ``top" and ``bottom" faces, forming a 2-sphere.
In fact, $\wh{C}_6^2$ is the boundary of an octahedron.
Since $\wh{C}_6^3$ is a 5-simplex and hence contractible, we have that $\PD_2(C_6)=\{(2,3)\}$.

A more general explanation of this 2-sphere topology in $\wh{C}_6^2$ is that $C_6^2$ contains all possible edges except that it is missing edges between the ``antipodal" vertices.
Therefore the clique complex $\wh{C}_6^2$ can be thought of as the boundary of the cross-polytope on six vertices in $\R^3$, i.e., the boundary of the convex hull of the six vertices
\[(\pm1,0,0),(0,\pm1,0),(0,0,\pm1).\]
Since this cross-polytope is a 3-ball, its boundary is homeomorphic to a 2-sphere.
\end{ex}

\begin{figure}[h]
\begin{center}
	$\begin{array}{c@{\hspace{.5in}}c}
    	\relabelbox  {\epsfysize=1.5in \epsfbox{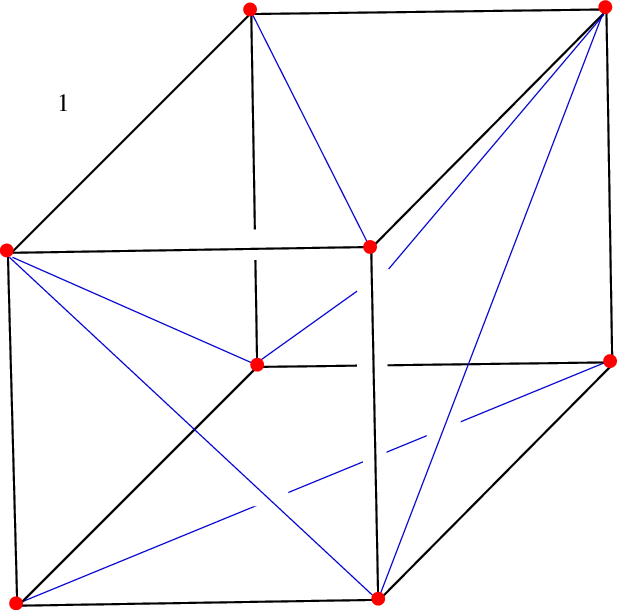}}
			\relabel{1}{\small $H$} 				
		\endrelabelbox &
		\relabelbox  {\epsfysize=1.5in \epsfbox{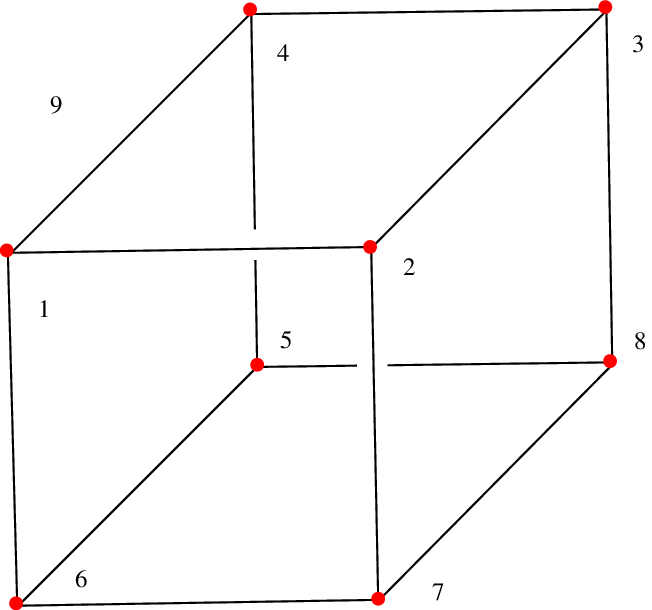}}
			\relabel{1}{\footnotesize $v_1$} 
			\relabel{2}{\footnotesize $v_2$}  
			\relabel{3}{\footnotesize $v_3$} 
			\relabel{4}{\footnotesize $v_4$} 
			\relabel{5}{\footnotesize $v_7$} 
			\relabel{6}{\footnotesize $v_6$} 
			\relabel{7}{\footnotesize $v_5$}  
			\relabel{8}{\footnotesize $v_8$} 
			\relabel{9}{\small $Q$} 				
		\endrelabelbox \\
	\end{array}$
\end{center}	
\caption{\footnotesize On the left, graph $H$ is the skeleton of a unit cube with diagonals in each face ($8$ vertices and $18$ edges).
On the right, the graph $Q$ is the 1-skeleton of a unit cube ($8$ vertices and $12$ edges).}
\label{fig-examples}
\end{figure}

\begin{ex}
Let $H$ be a graph with $8$ vertices and $18$ edges as shown in \cref{fig-examples}-left.
Then the simplicial complex $\wh{H}^1$, namely the clique complex of $H$, is topologically a sphere.
Since any two vertices are at distance at most $2$ apart in the shortest path metric on $G$, the clique complex $\wh{H}^2$ is a $7$-simplex, which is contractible.
This shows that $\PD_2(H)=\{(1,2)\}$.
Notice that $\PD_1(H)=\emptyset$ for this example.
\end{ex}

More generally than the prior example, let $\mathcal{T}$ be any simplicial complex triangulation of a $k$-sphere that is a clique complex, i.e., that is a flag triangulation of the sphere.
(For example, $\mathcal{T}$ cannot be the boundary of the $(k+1)$-simplex, as that is not a clique complex.)
If graph $G$ is the 1-skeleton of triangulation $\mathcal{T}$, then $\PD_k(G)$ contains a homology class born at $b=1$.

\begin{ex}
Let $Q$ be a graph with $8$ vertices and $12$ edges as shown in \cref{fig-examples}-right.
Graph $Q$ is the 1-skeleton of a unit cube.
Notice that $\wh{Q}^1=Q$ as there are no triangles in $Q$.

Interestingly, $\wh{Q}^2$ is homeomorphic to the 3-sphere $S^3$.
In order to see this, note that $Q^2$ contains all possible edges except that it is missing edges between the ``antipodal" vertices $v_1$ and $v_8$, $v_2$ and $v_7$, $v_3$ and $v_6$, and $v_4$ and $v_5$.
Therefore the clique complex $\wh{Q}^2$ is the boundary of the cross-polytope on 8 vertices in $\R^4$.
Since this cross-polytope is a 4-ball, its boundary is homeomorphic to a 3-sphere.
See~\cite{adamaszek2021vietoris,carlsson2020persistent,shukla2022vietoris} for further analyses of the clique complexes of hypercube graphs.

Since any two vertices are at distance at most $3$ apart in the shortest path metric on $Q$, it follows that $\wh{Q}^3$ is contractible.
We therefore have that $\PD_2(Q)$ is the empty diagram, while $\PD_3(Q)=\{(2,3)\}$.
\end{ex}

In the examples above, we have following persistence diagrams for 2-dimensional homology:
\[\PD_2(H)=\{(1,2)\}  \quad \quad \PD_2(C_6)=\{(2,3)\}.\]
These show that in higher dimensions ($k\ge 2$), the birth times may not be same for all $k$-cycles,  as we had in general for $\PD_1(G)$ with the 1-dimensional persistent homology of power filtrations of graphs.
So, a direct generalization of \cref{thm-main1} to higher homological dimensions is not possible.
While the examples above show that the birth times may not be constant for higher homologies, how large can the birth times be?

\begin{question} For $k\ge 2$, can $\PD_k(G)$ have an element with birth time $\ge k+1$ ?
\end{question}

The answer to this question is ``Yes.''
In~\cite[Corollary~6.7]{adamaszek2013clique}, Adamaszek gives a complete picture of the topological types of clique complexes of the powers of cyclic graphs.
Let $C_n$ be the cycle graph of length $n$.
Then, Adamaszek proves that the clique complex of the third power of $C_9$ has the homotopy type $\wh{C}_9^3\simeq \vee^2 S^2$, yielding a birth time $b=3$ in $\PD_2(C_9)$.
More generally, $\wh{C}_{3n}^n\simeq \vee^{n-1} S^2$ and $\wh{C}_{3n}^r\simeq S^1$ for $1\le r\le n-1$, which implies a birth time $b=n$ in $\PD_2(C_{3n})$.
In higher dimensions, $\wh{C}_{14}^5\simeq S^3$ which implies a birth time $b=5$ in $\PD_3(C_{14})$.
These examples show that the birth times for $k$-dimensional topological features can appear very late in the filtration for $k\geq 2$.
The constant birth time property is special to dimensions $k=0$ and $k=1$ for power filtrations of graphs.

These examples also show that an explicit description of $\PD_k(G)$, as we have when $k=0$ or when $k=1$ (see \cref{thm-main1}), is difficult for $k\ge2$.
However, in the following sections, we will discuss the next best thing: upper bounds for the life spans of higher-dimensional homology classes $\sigma$ of the form $d/b<C_\sigma$ for any feature $(b,d)\in \PD_k(G)$ corresponding to homology class $\sigma$, and we will give geometric interpretations of these bounds.

\begin{remark}[Multiplicative Persistence]
In persistent homology, the persistence (life span) of a $k$-homology class $\sigma$ with $(b, d)\in PD_k(\X)$ is typically defined as $d-b$, and it is interpreted as the size of the cavity that $\sigma$ represents in the point cloud (or given data) $\X$.
Intuitively, while the death time $d$ estimates the radius of the cavity (homology class), the birth time $b$ estimates how quickly the cavity forms, i.e.\ the closer the points are that form the homology class, the earlier the birth time.
There is, however, another notion called \textit{multiplicative persistence} defined as $d/b$ to describe the life spans of a topological feature $\sigma$~\cite{bobrowski2017maximally}.
This quantity is scale invariant, i.e.\ similar shapes in different scales will produce similar multiplicative persistence life spans.
For more discussion on multiplicative persistence, see~\cite[Section 3]{bobrowski2017maximally}.
In the literature, the use of multiplicative persistence (or persistence in logarithmic scale) is common~\cite{bobrowski2017maximally,chazal2013persistence, phillips2015geometric,sheehy2014persistent, buchet2016efficient}.
Multiplicative persistence is also better suited for the estimates in our main results (\cref{thm:main-width,thm:main-width-k}).
\end{remark}

\section{Conjectured Upper Bounds for Life Spans via Volume}
\label{sec:area}

In this section we describe a possible connection between thick-thin decompositions and life spans of persistent homology features.
We make conjectures regarding upper bounds on the death scales of persistent homology classes in terms of the volume of a minimal representative.
Most of the terms and ideas introduced in this section are well-established notions in geometric topology~\cite{thurston1997three,  colding2011course,benedetti2012lectures}.
They provide useful geometric intuition and the inspiration for the ideas in the next section, \cref{sec:width}, where we prove a bound by introducing a notion called the {\em width} of the cycle, which basically measures the ``thickness'' of the cycle.

Let $G$ be a graph, and let $\wh{G}^n$ be the clique complex of $G^n$, the $n^{th}$ power of $G$.
Let $\wh{G}^0\subset \wh{G}^1\subset \ldots \subset \wh{G}^D$ be the power filtration of $G$ as defined in \cref{ssec:powerfiltration}.
Let $\sigma$ be a $2$-dimensional homology class corresponding to the persistence diagram point $(b,d)\in\PD_2(G)$, and let
$S$ be a 2-cycle in $\wh{G}^b$ that generates $\sigma$ and that has as few $2$-simplices as possible.
We denote this number of $2$-simplices as $|\sigma|$, the \emph{area} of the homology class $\sigma$.
As an example, $S$ may be a genus-$g$ surface ($g\ge 0$).
But more generally $S$ need not be a manifold: $S$ could be a torus $S^1\times S^1$ with one component circle $S^1\times\{x\}$ collapsed to a point, or $S$ could be the wedge sum of two surfaces, etc.
We ask if one can provide upper bounds on the death scale $d$ by showing that $S$ is nullhomologous in $\wh{G}^n$, for some $n$ that depends on the geometry of $S$.

\begin{figure}[h]
\relabelbox  {\epsfxsize=4in
\centerline{\epsfbox{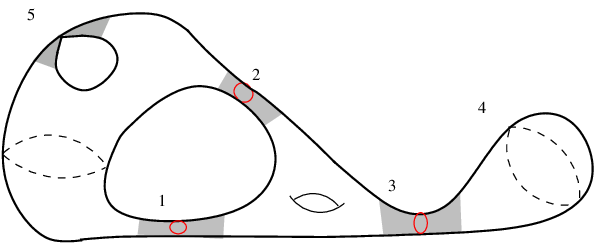}}}
\relabel{1}{\footnotesize $S^{thin}$} 
\relabel{2}{\footnotesize $S^{thin}$}  
\relabel{3}{\footnotesize $S^{thin}$} 
\relabel{4}{\small $S$} 
\relabel{5}{\footnotesize $S^{thin}$} 
\endrelabelbox
\caption{\textbf{Thick-Thin Decomposition.} $S$ has four thin components.}
\label{fig:thickthin}
\end{figure}

First, we describe a related geometric notion: a {\em thick-thin decomposition}, which decomposes a surface into its ``thick'' and ``thin'' parts; see \cref{fig:thickthin}.
The thick-thin decomposition is well-known in hyperbolic geometry in a different context, where it is called the Margulis Lemma~\cite{benedetti2012lectures}.
In that setting, the topology of the thin parts of a hyperbolic manifold is well-understood, and one can focus on the thick parts to understand the topology of the whole manifold.
We ask if a similar idea could be useful in our context.

\begin{figure}
\begin{center}
$\begin{array}{c@{\hspace{.4in}}c@{\hspace{.2in}}c}
	\relabelbox  {\epsfxsize=1.2in \epsfbox{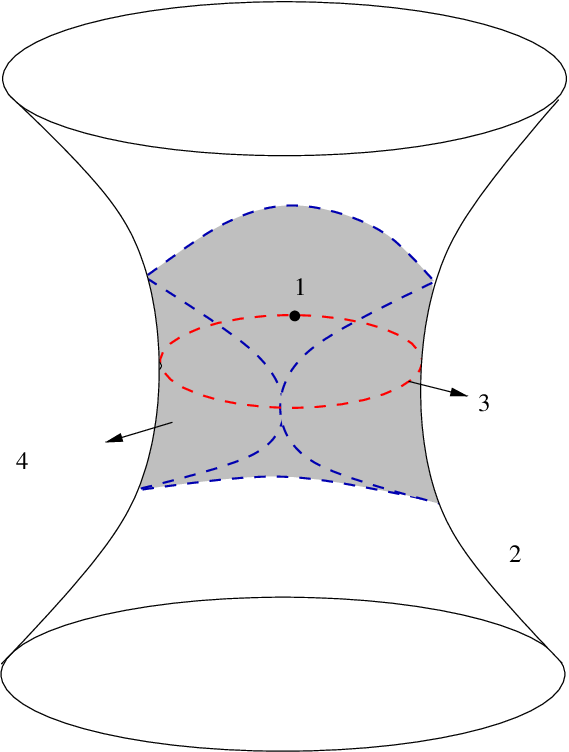}}
    \relabel{1}{\footnotesize $v$} 
    \relabel{2}{\footnotesize $S$}  
    \relabel{3}{\footnotesize bracelet} 
    \relabel{4}{\scriptsize $\mathrm{ball}(v)$} 
    \endrelabelbox
	&
	\relabelbox  {\epsfxsize=1.2in \epsfbox{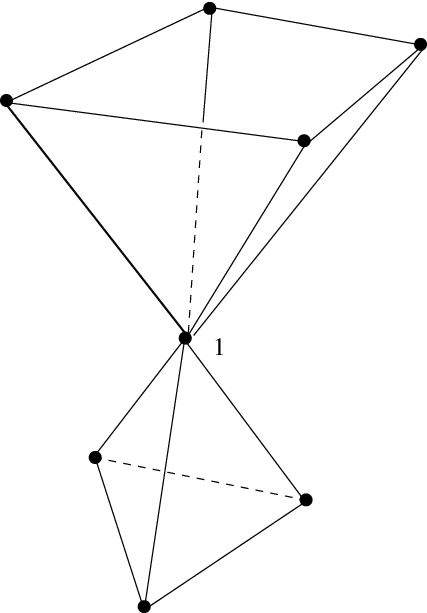}}
   	\relabel{1}{\footnotesize $v_0$} 
	\endrelabelbox
	%&
	%\relabelbox  {\epsfxsize=1.7in \epsfbox{figures2/bracelet2.eps}}
	%\relabel{1}{\footnotesize $v_0$} 
	%\relabel{2}{\scriptsize $w_1^-$}  
	%\relabel{3}{\scriptsize $w_2^-$} 
	%\relabel{4}{\scriptsize $w_3^-$} 
	%\relabel{5}{\scriptsize $w_4^-$} 
	%\relabel{6}{\scriptsize $w_5^-$}  
	%\relabel{7}{\scriptsize $w_5^+$} 	
	%\relabel{8}{\scriptsize $w_4^+$} 	
	%\relabel{9}{\scriptsize $w_3^+$} 	
	%\relabel{10}{\scriptsize $w_2^+$} 
	%\relabel{11}{\scriptsize $w_1^+$} 	
	%\endrelabelbox
\end{array}$
\end{center}
\caption{(Left) Ball about $v$ and a bracelet.
(Right) Pinch point.
}
\label{fig:inject1-pinch-bracelet}
\end{figure}

One can consider a thick-thin decomposition of our $2$-cycle $S$ as follows.
Roughly speaking, the injectivity radius $\inj(v)$ at a vertex $v\in S$ could be defined as the radius $r$ so that the ball in $S$ about $v$ radius $r$ does not ``overlap with itself'' in $S$, i.e.\ as the largest integer such that that the ball in $S$ about $v$ of radius $r$ is homeomorphic to a disc for all $1\le r\le \inj(v)$.
See \cref{fig:inject1-pinch-bracelet}-left.
We say that $w$ is a pinch point of $S$ if $\inj(w)=0$, as shown in \cref{fig:inject1-pinch-bracelet}-right.
We can now decompose $S$ into two parts, $S=\Sigma^-\cup\Sigma^+$, where the thin part $\Sigma^-$ contains all vertices with injectivity radius $\inj(v)\le \sqrt{|\sigma|}$, and where the thick part $S_m^{thick}$ contains all vertices with $\inj(v)\ge \sqrt{|\sigma|}$.

\begin{figure}[h]
\begin{center}
$\begin{array}{c@{\hspace{.2in}}c}
	\relabelbox  {\epsfxsize=2.2in \epsfbox{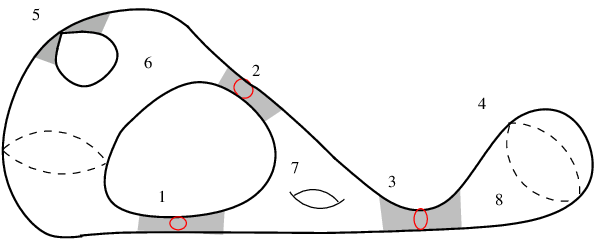}}
	\relabel{1}{\scriptsize $\Sigma_1^-$} 
	\relabel{2}{\scriptsize $\Sigma_2^-$}  
	\relabel{3}{\scriptsize $\Sigma_3^-$}
	\relabel{5}{\scriptsize $\Sigma_4^-$} 
	\relabel{4}{\small $S$} 
	\relabel{6}{\scriptsize $\Sigma^+_1$} 
	\relabel{7}{\scriptsize $\Sigma^+_2$}  
	\relabel{8}{\scriptsize $\Sigma^+_3$} 	
	\endrelabelbox
	&
	\relabelbox  {\epsfxsize=2.5in \epsfbox{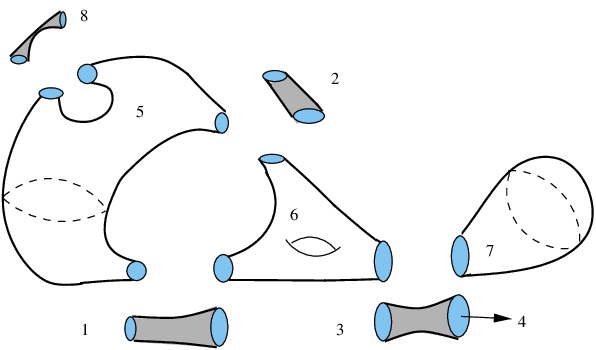}}
	\relabel{1}{\scriptsize $\wh{\Sigma}^-_1$} 
	\relabel{2}{\scriptsize $\wh{\Sigma}^-_2$}  
	\relabel{3}{\scriptsize $\wh{\Sigma}^-_3$}
	\relabel{8}{\scriptsize $\wh{\Sigma}^-_4$} 
	\relabel{4}{} 
	\relabel{5}{\scriptsize $\wh{\Sigma}^+_1$} 
	\relabel{6}{\scriptsize $\wh{\Sigma}^+_2$}  
	\relabel{7}{\scriptsize $\wh{\Sigma}^+_3$} 	
	\endrelabelbox
\end{array}$
\end{center}
\caption{Thick-Thin Decomposition of surface $S$.
Capping off each $\Sigma^\pm_i$ with blue disks in $\wh{G}^n$ produces each 2-cycle $\wh{\Sigma}^\pm_i$.}
\label{fig:thickthin2}
\end{figure}

Further decompose $\Sigma^-$ and $\Sigma^+$ into connected pieces $\Sigma^-_j$, and $\Sigma^+_k$.
Let $n\approx \sqrt{|\sigma|}\cdot b$.
We would like to cap off each piece to obtain the 2-cycles $\wh{\Sigma}_j^-$ and $\wh{\Sigma}^+_k$, which we would like to show are nullhomologous in $\wh{G}^n$.
Then, the refined decomposition $S=\sum_j\wh{\Sigma}^-_j +\sum_k\wh{\Sigma}^+_k$ as a sum of nullhomologous 2-cycles would give that $S$ is nullhomologous in $\wh{G}^n$.
See \cref{fig:thickthin2}.
Each ball of radius greater than $\sqrt{|\sigma|}$ about a vertex $v$ in a thin subsurface $\Sigma^-_j$ ``overlaps with itself'', creating at least one loop in this ball of bounded circumference (see \cref{fig:inject1-pinch-bracelet}-left).
A name for such a loop is a \emph{bracelet}, since it formed by wrapping a (topologically trivial) ball until it overlaps itself, in the same way that some bracelets are topological intervals wrapped around to form a circle on one's wrist.
We would like show that the vertex set of each such bracelet forms a complete simplex in $\wh{G}^n$, and that these simplices glue together to form a contractible simplicial complex, meaning that $\wh{\Sigma}^-_j$ is nullhomologous in $\wh{G}^{n}$.
Handling the pinch points will require care.
For the thick parts, since balls of radius at most $\sqrt{|\sigma|}$ about a vertex $v$ in a thick portion are embedded, one would like to prove an isoperimetric inequality lower bounding the number of 2-simplices in such a ball.
Then, once such a ball contains as many 2-simplices as there are in $\Sigma^+_k$, this provides an upper bound on the diameter of $\Sigma^+_k$, potentially showing that $\wh{\Sigma}^+_k$ will be nullhomologous in $\wh{G}^{n}$.
This leads us to the following question.

\begin{question}
\label{ques:area}
Let $(b,d)\in \PD_2(G)$ correspond to the 2-dimensional homology class $\sigma$.
Let $|\sigma|$ be the minimal number of 2-simplices in a 2-cycle generating this homology class.
Then is $d \le C_2 \sqrt{|\sigma|}\cdot b$ for some constant $C_2$?
%$d \le \left(\sqrt{|\sigma|}+1\right)b$?
\end{question}

Similarly, one could try to adapt the notion of thick-thin decompositions to higher dimensions, and try to use the volume to bound the diameter of equatorial hypersurfaces in the thick parts.

\begin{question}
\label{ques:vol2}
Let $(b,d)\in \PD_k(G)$ correspond to the $k$-dimensional homology class $\Omega$.
Let $|\Omega|$ be the volume of this homology class, i.e.\ the minimal number of $k$-simplices in a $k$-cycle generating this homology class.
Then is $d \le C_k \sqrt[k]{|\Omega|}\cdot b$ for some constant $C_k$ depending only on $k$?
%$d \le \left(\sqrt[k]{|\Omega|}+1\right)b$? 
\end{question}

We end this section with two simple examples of thick and thin surfaces.
These examples illustrate the relationship between the life span and the thickness of the surfaces.
One example shows that we can have large area surfaces with very short life spans in the power filtration, so area could only possibly be used as an upper bound for life span, and not as a lower bound.

\begin{figure}[h]
\centering
\begin{subfigure}{2.3in}
\centering
\includegraphics[width=2in]{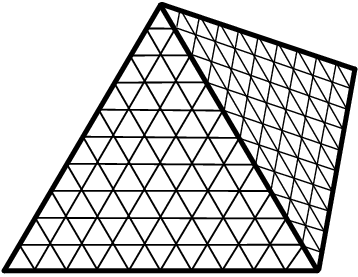}
\caption{Tetrahedron}
\end{subfigure}
\hspace{.4in}
\begin{subfigure}{2.3in}
\centering
		\relabelbox  {\epsfxsize=2.2in \epsfbox{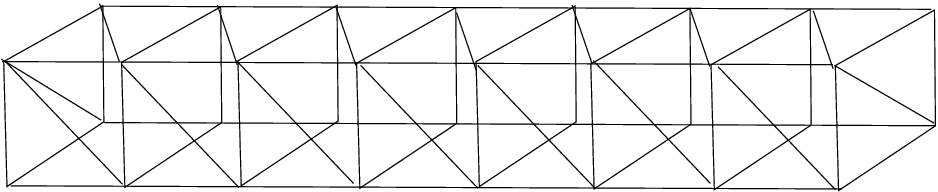}} \endrelabelbox
\vspace{.4cm}
\caption{Thin long box}
\end{subfigure}
\caption{\footnotesize \textbf{Thick and Thin Surfaces.}
The graph $G$ underlying the tetrahedron gives a thick surface, and has a long life span comparable to $\sqrt{\mbox{Area}}$.
The graph $G$ inducing a thin long box has a short life span in $\PD_2(G)$, even though it has very large area.}
\label{fig:boxes}
\end{figure}

\begin{ex}[Thick Surface and Long Life Span]
\label{ex:thick}
Consider the example in \cref{fig:boxes}-left: A tetrahedron $\sigma$ of edge length $m$.
This is an example of a thick surface whose injectivity radius is large, and depending on $m$.
The total area is $|\sigma|=4m^2$.
The $2$-cycle $\sigma$ is born at $b=1$, and its death scale is at least as large as $d\ge\frac{m}{2}$, meaning that its multiplicative life span is at least $d/b\ge\frac{m}{2}$.

We give a sketch (but not a complete proof) of why the death scale is at least as large as $d\ge\frac{m}{2}$.
Embed the graph $G$ in $\mathbb{R}^3$ as the 1-skeleton of a (triangulated) regular tetrahedron whose barycenter is at the origin in $\mathbb{R}^3$.
By extending linearly to simplices, this provides a map $\wh{G}^n\to \mathbb{R}^3$ for all $n\ge0$; this map $\wh{G}^n\to \mathbb{R}^3$ is not an embedding once $\wh{G}^n$ contains simplices above dimension $3$.
One can show that the death time $d$ is at least as large as the smallest scale parameter $n$ such that the map $\wh{G}^n\to \mathbb{R}^3$ hits the origin $\vec{0}$ in $\mathbb{R}^3$.
Indeed if this map misses the origin, then up to homotopy it is a retract onto $\mathbb{R}^3\setminus\{\vec{0}\}\simeq S^2$, and therefore is surjective onto the 2-dimensional homology of $S^2$.
(See~\cite[Proposition~5.3]{AAF} and~\cite[Corollary~21]{ABF} for related ideas.)
Furthermore, one can show that the smallest $n$ such that the map $\wh{G}^n\to \mathbb{R}^3$ hits the origin in $\mathbb{R}^3$ is the diameter of a (smaller) tetrahedron whose vertices are placed at or near the four centers of each of the four 2-dimensional faces of the (larger) tetrahedron represented by graph $G$, and that this smaller tetrahedron has diameter at least as large as $m/2$.
\end{ex}

\begin{ex}[Thin Surface and Short Life Span]
\label{ex:thin}
The second example is given in \cref{fig:boxes}-right: a thin long box  $\sigma$ of dimension $1\times 1\times m$.
This is an example of a thin surface with small injectivity radius for all vertices not at the ends of the cylinder.
The total area (number of triangles) of $\sigma$ is $8m+4$, i.e.\ $|\sigma|=8m+4$.
However, for any $m$ a $2$-cycle $\sigma$ is born at $b=1$, and dies at $d=2$, i.e.\ $(b,d)= (1,2)$.
Hence, its multiplicative life span is $d/b=2$.
However, in \cref{sec:width} we will describe a new geometric notion, the ``width'' of a surface, which will provide a better upper bound for the life span of this 2-cycle.
\end{ex}

The example of the thin long box shows that an upper bound on the life span based on area leaves room for improvement.
In the following section, we pursue an idea motivated by this example: upper bounding the life span of a persistent homology feature using the ``width'' of a generator.

\section{Upper Bound for Life Spans via Width}
\label{sec:width}

There is another natural way to define the thickness of a surface in geometric topology: the \textit{width}~\cite{gromov1988width}.
In the literature, this is also called the Urysohn width (see \cref{rmk:width}).
There are several ways to define this notion.
Here, we use sweepouts and the min-max technique to give an easy introduction to this concept.
By applying this idea to our context, we prove an upper bound on the life span of persistent homology features.
In this section we restrict attention to 2-dimensional homology, before considering $k$-dimensional homology in \cref{sec:higher}.
Throughout the paper, we use $\Z/2\Z$-coefficients for homology for simplicity, but our proofs can be adapted to the other field coefficients with the corresponding modifications.

First, we need to describe a metric that we use in our proof.
Let $G$ be a finite simple graph that is connected, and let $\V$ be its vertex set.
Recall that we have defined a metric $\rho\colon \V\times\V\to \R$ by letting $\rho(v,v')$ be the length of the shortest path in $G$ between $v$ and $v'$.
A related metric is defined by $\rho_n(v,v')=\lceil \tfrac{\rho(v,v')}{n}\rceil$.
To see that $\rho_n$ indeed defines a metric, note that $\rho_n(v,v')$ is simply the length of the shortest path between $v$ and $v'$ in the graph in $G^n$.
Clearly $\rho_1=\rho$, but typically we get a different metric $\rho_n$ whenever $n\ge 2$.
For example, let $G$ be the graph from \cref{fig-examples}-right, and let $n=2$.
Then using the vertex labels from \cref{fig-examples}-right, we have $\rho(v_1,v_3)=2$, whereas $\rho_2(v_1,v_3)=1$ since there is an edge from $v_1$ to $v_3$ in $G^2$.
The metric $\rho_n$ is the metric that will be most relevant in \cref{sec:width,sec:higher}.

From geometric point of view, the difference between the metrics $\rho$ and $\rho_n$ is related to the notions of \textit{intrinsic} and \textit{extrinsic} distances of a subspace inside a larger ambient space.
For example, consider the unit circle as a subspace of the larger ambient space $\R^2$.
In the circle, the intrinsic distance between a point and its antipode is $\pi$, which is the length of the path between these points in the circle.
By contrast, the extrinsic distance between two antipodal points is equal to $2$, the diameter of the circle in $\R^2$.
By definition, the intrinsic distance is always less than or equal to the extrinsic distance.
In our case, we have a nested sequence of increasing spaces $\wh{G}^1\subset \wh{G}^2\subset\ldots\subset \wh{G}^n\subset\ldots\subset \wh{G}^D$.
As $n$ increases, the space of possible paths between two vertices in the vertex set $\V$ increase, and therefore it may be possible to find shorter paths in the larger space.

\subsection{Min-Max Technique}

The min-max technique is basically a way to measure the thickness of a surface (or $k$-manifold).
The idea is to look at the surface from each direction, and to find the direction in which the surface looks ``the thinnest.''
In order to look at the surface from each direction, we define sweepouts~\cite{colding2003min}.

Let $\sigma$ be a $2$-dimensional homology class (with $\Z/2\Z$ coefficients) with birth and death times $(b,d)\in \PD_k(G)$.
Then, we can represent $\sigma$ by some $2$-dimensional cycle $S$ in $\wh{G}^{b}$ that generates the homology class $\sigma$.
Before defining the width of a homology class, we need to define the width of these cycles.

Let $\V(S)$ be the set of vertices that are in at least one $2$-simplex of $S$.
For $v\in \V(S)$, let $f_v\colon \V(S)\to \R$ be the function defined by $f_v(w)=\rho_b(v,w)$, where $\rho_b$ is the shortest path distance in $G^b$ (recall $\rho_b(v,w):=\lceil \tfrac{\rho(v,w)}{b}\rceil$).
Note $f_v^{-1}(i)$ is the set of all vertices in $\V(S)$ at distance $i$ (with respect to $\rho_b$) from $v$.
Let the eccentricity $e_v$ be the distance to a farthest vertex from $v$, namely $e_v=\max_{w\in \V(S)}\rho_b(v,w)$.

\begin{figure}[h]
	\relabelbox  {\epsfxsize=3.5in
		\centerline{\epsfbox{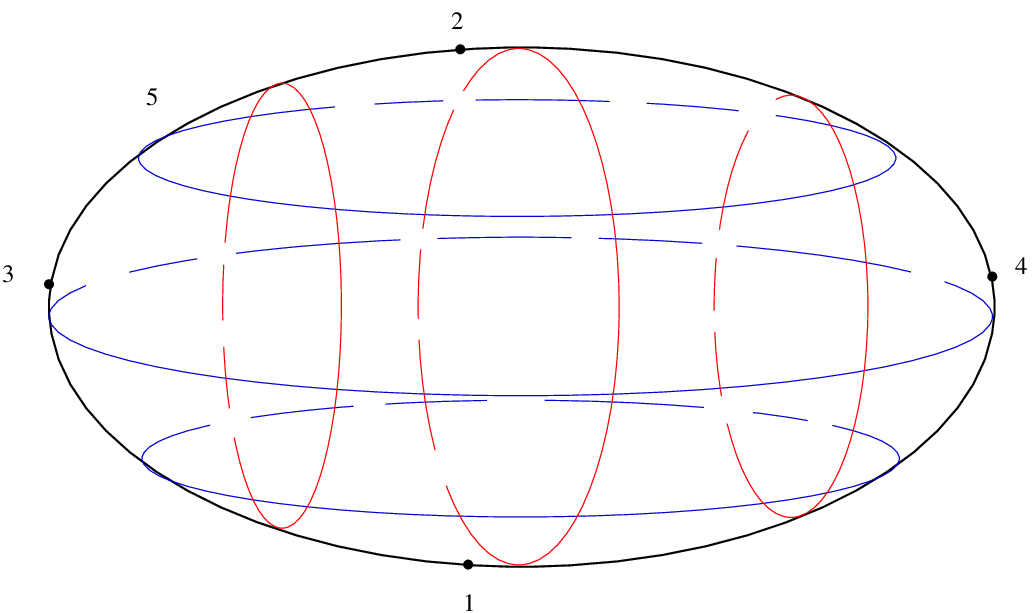}}}
		\relabel{1}{} 
		\relabel{2}{\footnotesize $v$}  
		\relabel{3}{\footnotesize $w$} 
		\relabel{4}{} 
		\relabel{5}{\footnotesize $S$} 
	\endrelabelbox
\caption{\textbf{Sweepouts.} Blue sweepout $\Lambda(v)=\{f_v^{-1}(i)\}_i$, and red sweepout $\Lambda(w)=\{f_w^{-1}(i)\}_i$.}
\label{fig:minmax}
\end{figure}

\begin{definition}[Sweepouts]
Let $S$ be a 2-cycle in $\wh{G}^b$ and let $v\in\V(S)$ be a vertex.
The \emph{sweepout of $S$}, $\Lambda(v)$, partitions $\V(S)$ into collections of vertices at distance $i$ from $v$, namely $\Lambda(v)=\{f_v^{-1}(i) \mid 0\le i\le e_v\}$.
\end{definition}

The collection of slices $f_v^{-1}(i)$ in a sweepout cuts $S$ into layers; see \cref{fig:minmax}.
Readers familiar with applied topology will recognize some similarities with sublevelset persistent homology, although with sweepouts the focus is on levelsets moreso than on sublevelsets.

For $Y\subseteq\V$, we define the \emph{diameter} of $Y$ to be the distance between the farthest two vertices in $Y$, namely $\diam(Y):=\max_{v,v'\in Y}\rho_b(v,v')$.
Note that if $\diam(Y) \le t$ for some $t>0$, then $\rho_b(v,v')\le t$ for all $v,v'\in Y$, meaning that $\rho(v,v')\le t\cdot b$ for all $v,v'\in Y$.
This means that the finite set $Y$ is simplex in $\wh{G}^{t\cdot b}$.
We will use this fact in the proof of \cref{thm:main-width}.

By using the sweepouts of $S$, we define the width of $S$~\cite{colding2005estimates}, and hence the width of a homology class, as follows:

\begin{definition}[Width of a Homology Class]
Let $S$ be a 2-cycle in $\wh{G}^b$.
Define the \emph{width} $\omega(S)$ as the minimum, over all sweepouts of $S$, of the maximal diameter of a slice in the sweepout:
\[\omega(S)=\min_{v\in\V(S)} \max_{\Lambda(v)} \{\diam(f_v^{-1}(i))\}.\]
Let $\sigma$ be a homology class in $H_2(\wh{G}^b)$.
Then, we define the \emph{width} of $\sigma$ as \[\omega(\sigma)=\min_{S\in\sigma}\omega(S),\]
where the minimum is taken over all 2-cycles $S$ in $\wh{G}^b$ that generate the homology class $\sigma$.
\end{definition}

\begin{figure}[h]
	\relabelbox  {\epsfxsize=3.3in
		\centerline{\epsfbox{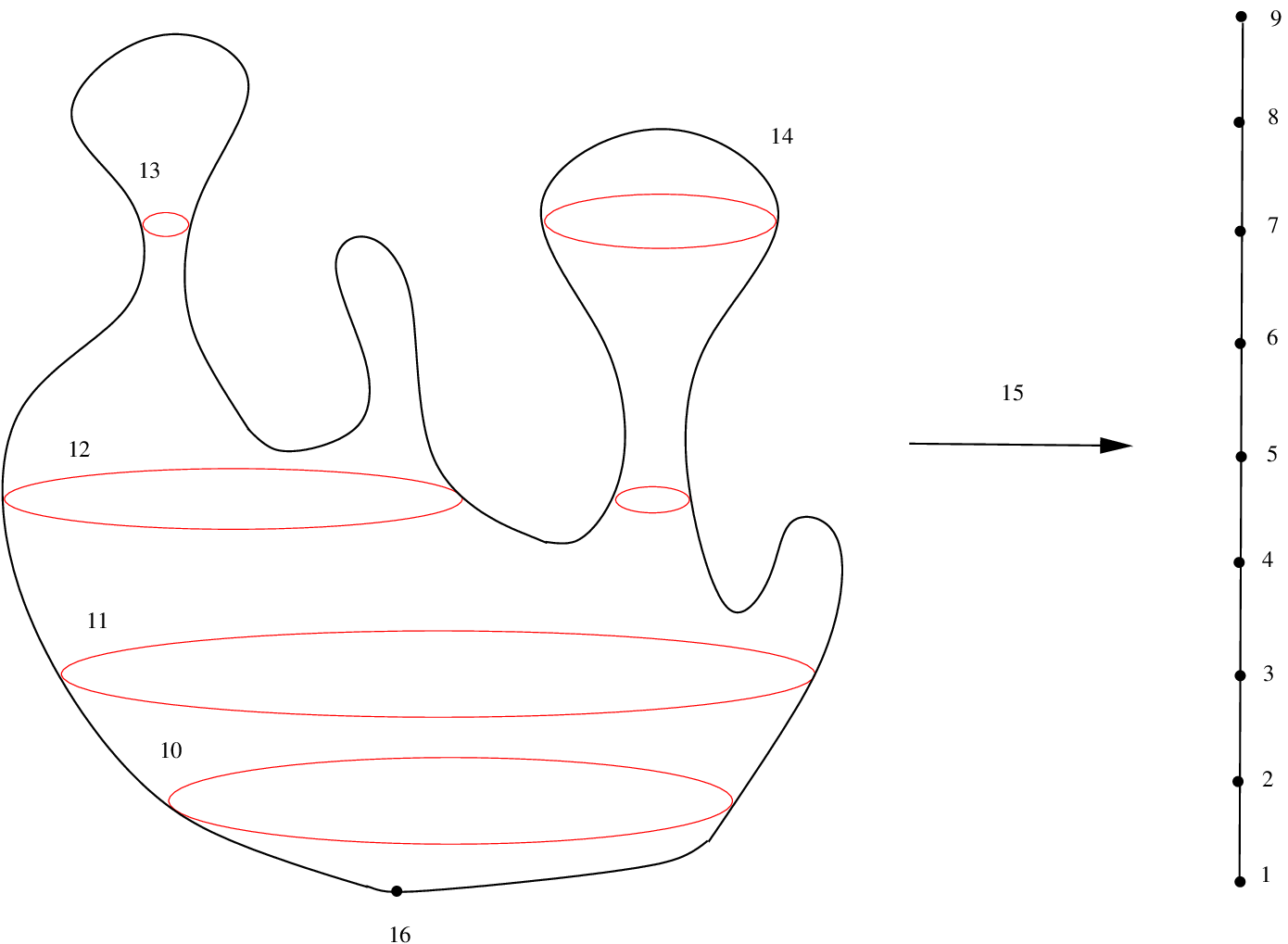}}}
		\relabel{1}{\footnotesize $0$} 
		\relabel{2}{\footnotesize $1$}  
		\relabel{3}{\footnotesize $2$} 
		\relabel{4}{\footnotesize $3$} 
		\relabel{5}{\footnotesize $4$} 
		\relabel{6}{\footnotesize $5$} 
		\relabel{7}{\footnotesize $6$}  
		\relabel{8}{\footnotesize $7$} 
		\relabel{9}{\footnotesize $8$} 
		\relabel{10}{\footnotesize $N_1$} 
		\relabel{11}{\footnotesize $N_2$}  
		\relabel{12}{\footnotesize $N_4$} 		
		\relabel{13}{\footnotesize $N_6$} 		
		\relabel{14}{\footnotesize $S$} 		
		\relabel{15}{\footnotesize $f_v$} 		
		\relabel{16}{\footnotesize $v$} 		
	\endrelabelbox
\caption{\textbf{Filter Function.} Slicing via a filter function: $N_i:=f_v^{-1}(i)$}
\label{fig:filtering}
\end{figure}

We are minimizing among maximum slices of $S$ in all sweepouts; see \cref{fig:minmax}.
This is called the min-max technique in minimal surface theory \cite{colding2003min}.
In metric geometry, this notion corresponds to the Urysohn $1$-width (see \cref{rmk:width}).

\begin{remark}[Subtleties in the Width Definition]
An important point in the width definition is that the diameter is defined using shortest paths in $G^b$, not using shortest paths in $G$, and not using shortest paths in the individual slices (which need not even be connected; see \cref{fig:filtering}).
\end{remark}

\subsection{Bounding Life Spans via Width}

The following theorem bounds the life span of a $2$-dimensional topological feature in terms of its width.
The main idea of the proof is to slice the manifold up into small diameter pieces (slabs) by using the best sweepout.

\begin{theorem}
\label{thm:main-width}
Let $(b,d)\in \PD_2(G)$ correspond to the $2$-dimensional homology class $\sigma$.
Let $\omega(\sigma)$ be the width of the 2-dimensional homology class $\sigma$ in $\wh{G}^b$.
Then, 
\[d \le (\omega(\sigma)+1) b.\]
\end{theorem}

%%%%%\note{Try to add an extra $\frac{1}{3}$ constant in the case of 2-dimensional homology?}
%%%\note{The notation in the commented-out paragraph below in the LaTeX source may be useful for adding the 1/3 constant.}
%Now, consider the graph $N_i$.
%Let $C^1_i, C^2_i,\ldots, C^{n_i}_i$ be loops which generate the 1-dimensional homology of $N_i$.
%Note that some loops might intersect each other.
%For each loop $C^j_i$, we define a least area disk $D^j_i$ as follows.
%Let $\{w_1,w_2,\ldots,w_{\lambda_i^j}\}$ be the vertices of $C^j_i$, where $|N_i^j|=\lambda_i^j$.
%Consider triangles $U_k=[w_1,w_k,w_{k+1}]$.
%Then, $D_i^j=\bigcup_{k=2}^{\lambda_i^j-1}U_k$ is a disk in $\wh{G}^{b+\omega(\sigma)+1}$ with $\partial D_i^j=N_i^j$.
%The reason $D_i^j$ exists in $\wh{G}^{b+\omega(\sigma)+1}$ is that for any $i$, $\diam(N_i)\le \omega(\sigma)$ by assumption.

\begin{proof}
Let $S$ be a 2-cycle in $\wh{G}^b$ minimizing the width of $\sigma$, i.e.\ $\omega(S)=\omega(\sigma)$.
Let $v\in\V(S)$ be such that $\Lambda(v)$ is a minimizing sweepout, meaning that we have $\diam(f_v^{-1}(i)) \le \omega(S)$ for all $0\le i\le e_v$, where $e_v$ is the eccentricity of $v$ in $S$.
We use the notation $N_i:=f_v^{-1}(i)$.
See \cref{fig:width}-left.

\begin{figure}[t]
\begin{center}
	$\begin{array}{c@{\hspace{.5in}}c}
	\relabelbox  {\epsfysize=2.8in \epsfbox{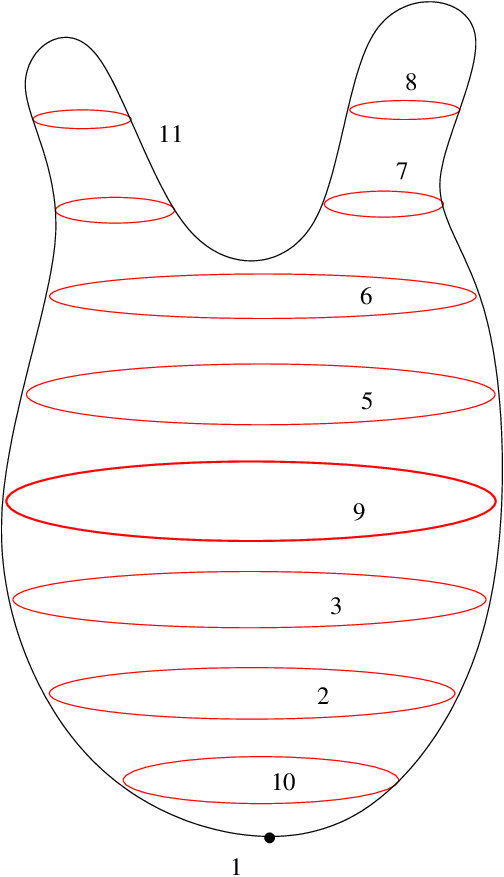}}
		\relabel{1}{\footnotesize $v_0$} 
		\relabel{10}{\footnotesize $N_1$} 
		\relabel{2}{\footnotesize $N_2$}  
		\relabel{3}{\footnotesize $N_3$} 
		\relabel{5}{\footnotesize $N_5$} 
		\relabel{6}{\footnotesize $N_6$} 
		\relabel{7}{\footnotesize $N_7$}  
		\relabel{8}{\footnotesize $N_8$} 
		\relabel{9}{\footnotesize $N_4$} 
		\relabel{11}{ $S$} 
	\endrelabelbox &
	\relabelbox  {\epsfysize=2.8in \epsfbox{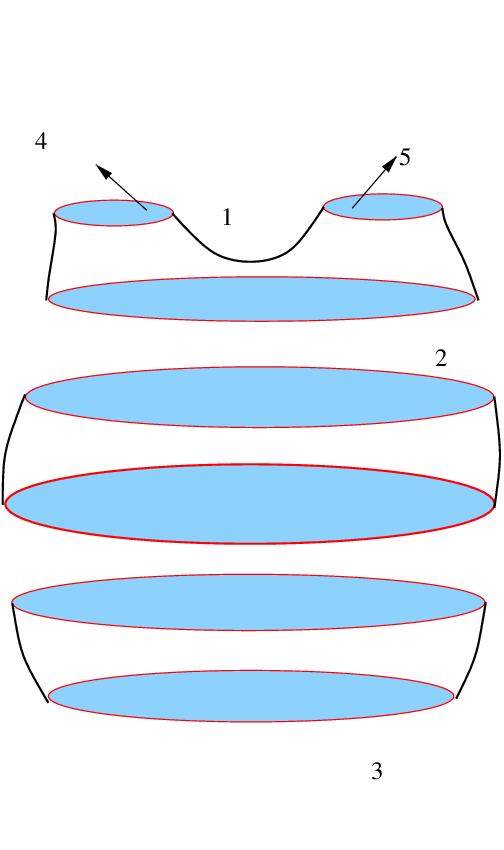}}
		\relabel{1}{ $\wh{\A}_7$} 
		\relabel{2}{ $\wh{\A}_5$}  
		\relabel{3}{ $\wh{\A}_3$} 
		\relabel{4}{ $\U_7$}  
		\relabel{5}{ $\U_7$} 
	\endrelabelbox \\
	\end{array}$
\end{center}
\caption{On the left, $\Lambda(v)$ is the sweepout with the smallest maximal slice $N_4$, i.e.\ $\omega(S)=\diam(N_4)$.
On the right, we cap off each slab $\A_i$ with the caps $\U_{i-1}$ and $\U_i$, to obtain the $2$-cycle $\wh{\A}_i$ in $\wh{G}^{(\omega(\sigma)+1)b}$.
As the $1$-cycle $Z_7$ associated to the slice $N_7$ contains two loops, note the $2$-chain $\U_7$ contains two disks.}
\label{fig:width}
\end{figure}

Let the \emph{slab} $\A_i$, for $1\le i\le e_v$, be the set of $2$-simplices in $S$ with vertices in $N_{i-1}\cup N_i$.
In particular, the slab $\A_i$ is the portion of $S$ between the \emph{slices} $N_{i-1}$ and $N_i$; see \cref{fig:width}.

We claim $S=\sum_{i=1}^{e_v}\A_i$ as a sum of 2-chains.
Indeed, any $2$-simplex $\tau$ in $S$ is also a simplex in $\wh{G}^b$, and therefore its vertices are at distance at most $1$ apart in the shortest-path metric $\rho_b$ on the graph $G^b$.
This implies that $\tau$ is in $\A_i$ for some $i$.
Indeed, for any two vertices $u,u'\in \tau$, we have $|f_v(u)-f_v(u')|=|\rho_b(v,u)-\rho_b(v,u')|\le \rho_b(u,u')\le 1$ by the triangle inequality.
The equality $S=\sum_{i=1}^{e_v}\A_i$ implies that $\emptyset = \partial S = \sum_{i=1}^{e_v} \partial\A_i$, where the sum is taken with coefficients in $\Z/2\Z$.

An edge $\tau$ can satisfy $\tau\in\partial\A_i\cap\partial\A_j$ for $i<j$ only if $j=i+1$, in which case $\tau$ is a subset of $N_i$.
Therefore the equality $\emptyset = \sum_{i=1}^{e_v} \partial\A_i$ implies that there exist $1$-cycles $Z_i$ for $0\le i\le e_v$ (with $Z_0=\emptyset=Z_{e_v}$) such that each edge in $Z_i$ is a subset of $N_i$, and $\partial\A_i=Z_{i-1} + Z_i$ for all $1\le i\le e_v$.

As $\diam(N_i)\le \omega(\sigma)$, $N_i$ generates a complete simplex in $\wh{G}^{\omega(\sigma)b} \subseteq \wh{G}^{(\omega(\sigma)+1)b}$.
As this simplex is contractible, the $1$-cycle $Z_i$ is the boundary of some $2$-chain $\U_i$ whose simplices are all subsets of $N_i$.

Define the $2$-cycles $\wh{\A}_i=\A_i\cup \U_{i-1}\cup \U_i$ for $1\le i\le e_v$.
The 2-chains $\U_{i-1}$ and $\U_i$ cap off the lower boundary $Z_{i-1}$ and the upper boundary $Z_i$ of $\A_i$; see \cref{fig:width}-right.
Note that $\wh{\A}_i$ is indeed a cycle, namely $\partial \wh{\A}_i=\emptyset$, which follows from the fact that $\partial\A_i=Z_{i-1}+Z_i=\partial\U_{i-1}+\partial\U_i$.
We have $S=\sum_{i=1}^{e_v}\wh{\A}_i$, since the caps $\U_i$ cancel in pairs.

We consider the diameter of $N_{i-1}\cup N_i$.
For any pair of vertices $u,w\in N_{i-1}\cup N_i$, either both belong to slice $N_{i-1}$, or both belong to slice $N_i$, or one vertex belongs to each.
If both vertices belong to the same $N_j$, then $\rho_b(u,w)\le \diam(N_j) \le \omega(\sigma)$.
If $u\in N_{i-1}$ and $w\in N_i$, then let $w'$ be a vertex in $N_{i-1}$ with $\rho_b(w,w')=1$.
Since $\rho_b(u,w')\le \omega(\sigma)$, we have $\rho_b(u,w)\le \omega(\sigma)+1$.
Therefore, $\diam(N_{i-1}\cup N_i)\le\omega(\sigma)+1$.

As $\diam(N_{i-1}\cup N_i)\le\omega(\Omega)+1$, the vertex set $N_{i-1}\cup N_i$ forms a simplex in the simplicial complex $\wh{G}^{(\omega(\sigma)+1)b}$.
Therefore, $\wh{\A}_i$ is nullhomologous in $\wh{G}^{(\omega(\sigma)+1)b}$.
Since $S=\sum_{i=1}^{e_v}\wh{\A}_i$ as 2-cycles, it follows that $S$ is nullhomologous in $\wh{G}^{(\omega(\sigma)+1)b}$.
So $d \le (\omega(\sigma)+1)b$.
\end{proof}

\begin{remark}[Slicing Technique]
Notice that the main idea of the proof above is to slice the surface $S$ into small diameter subsurface ``slabs'' $\{\A_i\}$.
We call this the ``Slicing Technique''.
We note that this is similar to ideas also used by Virk in the context of Vietoris--Rips complexes to prove~\cite[Theorem~7.1]{virk2021footprints}; see in particular Figure~6 within.
In the following section, we generalize the slicing technique to any dimension, and we use the width to prove an upper bound on the persistent homology life spans in general homological dimensions.
\end{remark}

\begin{remark}[Improving the Upper Bound]
\label{rmk:improve}
Notice that in both \cref{ques:area} and \cref{thm:main-width}, our main technical approach is to chop the original surface $S$ in $\wh{G}^b$ into small diameter pieces $\Sigma_i$ with $\diam(\Sigma_i)\le t$ for some $t>0$.
Then, as each piece generates a complete simplex on its vertices in $\wh{G}^{t\cdot b}$, it will be trivially nullhomologous in $\wh{G}^{t\cdot b}$.
Here, we use a rough estimate to get nullhomologous surfaces.
On the other hand, in~\cite[Proposition~9.1]{lim2020vietoris}, the authors prove that the ``spread'' gives an upper bound on the persistence $d-b$ for any persistence bar in the Vietoris--Rips filtration, where for a geodesic space $X$ one has $\mathrm{spread}(X)\leq \frac{2}{3}\diam(X)$. 
It might be possible to adapt this spread result for geodesic spaces to our discrete setting, and improve our estimates using a notion of width that considers the spread (instead of the diameter) of each slice.
\end{remark}

\section{$\PD_k(G)$: Generalization to Higher Dimensions}
\label{sec:higher}

The ideas in the proof \cref{thm:main-width} are suitable to generalize to higher dimensions, which we do in this section.

Let $\Omega$ be a $k$-dimensional homology class, with $\Z/2\Z$ coefficients, with birth and death times $(b,d)\in \PD_k(G)$.
Then, we can represent $\Omega$ by some $k$-dimensional cycle $M$ in $\wh{G}^b$ that generates the homology class $\Omega$.
Before defining the width of a homology class, we need to define the width of these cycles.

Let $\V(M)$ be the set of vertices contained in at least one $k$-simplex of $M$, and let $v\in\V(M)$.
Let $f_v\colon \V(M)\to \R$ be the function defined by $f_v(w)=\rho_b(v,w)$, where $\rho_b$ is the shortest path distance in the graph $G^b$.
Note $f_v^{-1}(i)$ is the set of all vertices in $\V(M)$ at distance $i$ from $v$.
Let the eccentricity $e_v$ be the distance to a farthest vertex from $v$, namely $e_v=\max_{w\in \V(M)}\rho_b(v,w)$.

The sweepout $\Lambda(v)$ partitions $\V(M)$ into collections of vertices at distance $i$ from $v$, namely $\Lambda(v)=\{f_v^{-1}(i) \mid 0\le i\le e_v\}$.
We define the \emph{width} of the $k$-dimensional cycle $M$ in $\wh{G}^b$ as the minimum, over all sweepouts, of the maximum diameter (using the metric $\rho_b$) of a slice in the sweepout: \[\omega(M)=\min_{v}\max_{\Lambda(v)}\{\diam(f_v^{-1}(i))\}.\]
As before, in order to define the width of a homology class, we minimize the width among the cycle representatives $M$ of $\Omega$ in $\wh{G}^b$.
That is,
$\omega(\Omega)=\min_{M\in\Omega}\omega(M)$.
We can give the generalization of the theorem as before.

\begin{theorem}
\label{thm:main-width-k}
Let $(b,d)\in \PD_k(G)$.
Let $\omega(\Omega)$ be the width of the $k$-\linebreak dimensional homology class $\Omega$ in $\wh{G}^b$.
Then, 
\[d \le (\omega(\Omega)+1) b.\]
\end{theorem}

\begin{proof}
We apply the ``Slicing Technique" developed in the proof of \cref{thm:main-width} to higher dimensions.
Let $M$ be a $k$-cycle minimizing the width of $\Omega$,
and let $v\in\V(M)$ be such that $\Lambda(v)$ is a minimizing sweepout of $M$.
We use the notation $N_i:=f_v^{-1}(i)$ for $0\le i\le e_v$.

For $1\le i\le e_v$, let the slab $\A_i$ be the set of $k$-simplices in $M$ with vertices in $N_{i-1}\cup N_i$.
To see that $M=\sum_{i=1}^{e_v}\A_i$, note that any $k$-simplex $\tau$ in $M$ is also a simplex in $\wh{G}^{b}$, and therefore its vertices are at distance at most $1$ apart in the shortest-path metric $\rho_b$.
By the triangle inequality, this implies that $\tau$ is in $\A_i$ for some $i$.
The equality $M=\sum_{i=1}^{e_v}\A_i$ implies that $\emptyset = \partial M = \sum_{i=1}^{e_v} \partial\A_i$, where the sum is taken with coefficients in $\Z/2\Z$.

A $(k-1)$-simplex $\tau$ can satisfy $\tau\in\partial\A_i\cap\partial\A_j$ for $i<j$ only if $j=i+1$, in which case $\tau$ is a subset of $N_i$.
Therefore the equality $\emptyset = \sum_{i=1}^{e_v} \partial\A_i$ implies that there exist $(k-1)$-cycles $Z_i$ for $0\le i\le e_v$ (with $Z_0=\emptyset=Z_{e_v}$) such that each simplex in $Z_i$ is a subset of $N_i$, and such that $\partial\A_i=Z_{i-1} + Z_i$ for all $1\le i\le e_v$. 

As $\diam(N_i)\le \omega(\Omega)$, $N_i$ generates a complete simplex in $\wh{G}^{\omega(\Omega)b} \subseteq \wh{G}^{(\omega(\Omega)+1)b}$.
Since this simplex is contractible, the $(k-1)$-cycle $Z_i$ is the boundary of some $k$-chain $\U_i$ whose simplices are all subsets of $N_i$.
 
We define the $k$-cycle $\wh{\A}_i=\A_i\cup \U_{i-1}\cup \U_i$, which caps off the lower boundary $Z_{i-1}$ and the upper boundary $Z_i$ of $\A_i$.
Note $\wh{\A}_i$ is indeed a cycle since $\partial\A_i=Z_{i-1}+Z_i=\partial\U_{i-1}+\partial\U_i$ implies $\partial\wh{\A}_i=\emptyset$.
We have $M=\sum_{i=1}^{e_v}\wh{\A}_i$, since the caps $\U_i$ cancel in pairs.

We claim that $\diam(N_{i-1}\cup N_i)\le\omega(\Omega)+1$.
For any pair of vertices $u,w\in N_{i-1}\cup N_i$, either both belong to slice $N_{i-1}$, or both belong to slice $N_i$, or one vertex belongs to each.
If both vertices belong to the same $N_j$, then $\rho_b(u,w)\le \diam(N_j) \le \omega(\Omega)$.
If $u\in N_{i-1}$ and $w\in N_i$, then let $w'$ be a vertex in $N_{i-1}$ with $\rho_b(w,w')=1$.
Since $\rho_b(u,w')\le \omega(\Omega)$, we have $\rho_b(u,w)\le \omega(\Omega)+1$.

As $\diam(N_{i-1}\cup N_i)\le\omega(\Omega)+1$, the vertex set $N_{i-1}\cup N_i$ forms a simplex in the simplicial complex $\wh{G}^{(\omega(\Omega)+1)b}$.
Therefore, each $\wh{\A}_i$ is nullhomologous in $\wh{G}^{(\omega(\Omega)+1)b}$, and so $M=\sum_{i=1}^{e_v}\wh{\A}_i$ is nullhomologous in $\wh{G}^{(\omega(\Omega)+1)b}$.
It follows that $d \le (\omega(\Omega)+1)b$.
\end{proof}

\begin{remark}[Urysohn Widths]
\label{rmk:width}
Notice that in the above theorem, we used the same width notion defined for 2-cycles, now in the context of $k$-cycles.
For a $k$-manifold $M$, this is called Urysohn $1$-width $\omega_1(M)$, as it is induced by the codimension-$1$ sweepouts $f^{-1}(t)$ induced by $f\colon M\to\R$.
If one use instead a function $F\colon M\to\R^m$, then the sweepouts $F^{-1}(\mathbf{x})$ will be codimension-$m$ submanifolds; a similar definition gives us the Urysohn $m$-width $\omega_m(M)$ of the manifold~\cite{gromov1988width}.
It is known that for a closed $k$-manifold, the widths are monotone, namely
$\omega_1(M)\ge \omega_2(M)\ge\ldots\ge\omega_{k-1}(M)\ge \omega_k(M)=0$.
It is conjectured that the $(k-1)$-width $\omega_{k-1}(M)$ is closely related to Gromov's filling radius~\cite{guth2011volumes}.
So, by generalizing the arguments above, one might get better estimates by using $\omega_{k-1}(M)$.
See \cref{ssec:lowerbound} for further discussion.
\end{remark}

\section{Final Remarks}
\label{sec:remarks}

Our aim in this paper is two-fold.
First, we aim to bring fresh ideas from geometric topology to address subtle questions in applied topology and topological data analysis.
Second, we hope to attract the attention of geometric topologists to the emerging field of applied topology by showing how these tools are effective, as there are many more problems to be tackled.
In the following, we give some concluding remarks, and discuss further directions.

\subsection{Geometric Interpretation of $\PD_k(G)$}
\label{ssec:interpret}

In this paper, for a given $k$-dimensional topological feature (homology class) $\Omega$ in the power filtration of a graph $G$, we study the relation between its representation $(b,d)$ in the persistence diagram $\PD_k(G)$ and its geometric size.
Here, the geometric size can be interpreted as the volume or as the width of the homology class $\Omega$.

Upon going over the ideas in \cref{sec:area,sec:width,sec:higher}, one sees that a persistent homology class can persist as much as its size allows it to.
Even though in dimension $1$ the length of the homology class gives a good interpretation of this notion of size, in higher dimensions the volume can be a weak way to measure this.
A better notion is the width, which one can consider as a \emph{filling diameter} of the homology class.
The thick and thin surfaces in \cref{ex:thick} and \cref{ex:thin} are the key to understanding the width idea as the filling diameter of a homology class.
In the thin long box example, even though the area is large, the life span is very small.
The reason is that the width better measures the size of the $3$-dimensional body filling this surface.
In the thin long box case, the diameter of a largest embedded ball in this filling $3$-dimensional body is very small, which informally represents the width idea we are describing here.
In particular, since we do not have a well-defined way to describe this filling $3$-body, we define the width by using the geometry of the surface.
From this point of view, the width we use can informally be related to Gromov's notion of the filling radius \cite{gromov1983filling}, also used in \cite{lim2020vietoris}.

\subsection{Lower Bounds for Life Spans}
\label{ssec:lowerbound}

In this paper, we only give upper bounds for the life spans of homology classes in persistent homology.
Of course, the next natural question is the lower bounds via area or width, or via some new geometric notion.

Unfortunately, as the thin long box example in \cref{ex:thick} shows, the area cannot be used for a lower bound for the life spans.
One can make the area of the thin parts of the surface as big as one wants without effecting its life span.

However, the width would be a good candidate to get a lower bound for the life span.
As we discussed in previous section, the life spans are directly related to the notion of the filling diameter of the surface (or $k$-dimensional homology class).
It is believed that the width and Gromov's filling radius of a closed manifold are closely related~\cite{guth2011volumes, gromov1988width} (see \cref{rmk:width}).
In particular, the life span informally represents the diameter of the largest embedded ball in the $3$-dimensional body filling our surface.
From that perspective, the width notion mostly captures this idea, potentially giving route towards a lower bound via width.
In our paper, we use width in a simple way to chop the surface into small diameter surfaces.
However, it might be possible to adapt higher Urysohn widths notion to this context (\cref{rmk:width}) to get a lower bound for life spans, by injecting ideas from metric geometry where higher Urysohn widths are used to obtain lower bound for Gromov's filling radius~\cite{nabutovsky2021sweepouts}.

\subsection{Generalization to Point Clouds and Other Filtrations}

In the language of Vietoris--Rips complexes, the main result of our paper can be summarized as follows.
If $G$ is an unweighted connected finite graph, then in the persistent homology of the Vietoris--Rips complex of its vertex set, the death time $d$ of a $k$-dimensional homology class $\Omega$ is at most $(\omega(\Omega)+1) b$, where $\omega(\Omega)$ is the width of $\Omega$ (\cref{thm:main-width-k}).
As one can notice, many geometric ideas introduced in this paper, like the thick-thin decomposition, min-max, or sweepouts, can be generalized to the point cloud setting, or to other filtrations.
We therefore expect there to be generalizations of the above result to the setting of Vietoris--Rips, \v{C}ech, or witness complexes of more general point clouds.
Depending on the filtration, the interpretation and measure of size of these quantities bounding the life span of a homology class will be different, but nevertheless we expect that the notions developed in this paper will help to effectively summarize the information obtained in persistence diagrams.

In order to generalize these techniques to the point cloud setting, the first main obstacle to overcome is the varying distances between points in the cloud.
In our case, we took a simple discrete metric for our point cloud (vertices of the graph) induced by shortest paths in an unweighted graph.
However, recent advances in coarse geometry could be the key to generalize these geometric notions to more general metric spaces and point clouds~\cite{burago2001course, burago2015uniform}.
For example, suppose that $X$ is a geodesic metric space that can be closely approximated in the Gromov-Hausdorff distance by a graph $G$ in which all edges have the same length.
Then the stability of persistent homology~\cite{ChazalDeSilvaOudot2014,chazal2009gromov} implies that the persistent homology of the Vietoris--Rips complexes of $X$ are close to the persistent homology of the Vietoris--Rips complexes of $G$, which is governed by our results.

The paper~\cite{lim2020vietoris} connects the Vietoris--Rips filtration to the filling radius of a space.
Given a metric space $X$, this same paper proves that any persistent homology bar in the Vietoris--Rips filtration of a metric space $X$ has a length bounded from above by the spread of $X$.
As the spread only depends on the metric space $X$, it gives the same bound for all homology classes (in all homological dimensions).
In our work, in the setting of unweighted graphs, we instead give refined bounds depending on the size of each individual homology class.
We speculate that some of our techniques may prove useful in the setting of more general metric spaces.

\subsection{TDA on Graphs} 
\label{ssec:tdaongraphs}

We focused on the setting of power filtrations of graphs to first study persistent homology in a discrete setting, in order to learn what might happen in more general contexts.
That said, our results show that higher persistence diagrams carry valuable information about graph properties.
In particular, even elementary graphs like the cyclic graphs~\cite{adamaszek2013clique}, or the small examples in \cref{sec:PD2}, have nontrivial higher-dimensional persistence diagrams.
These higher persistence diagrams can be effective to detect hidden patterns in real life datasets, such as widths of $k$-dimensional cavities.

Because of the computational cost of the power filtration, researchers applying TDA on graphs sometimes prefer sublevel (superlevel) filtrations defined by a filter function to study graph properties.
Instead, the power filtration considers the graph vertices as a point cloud where the distances are given by the graph, meaning that the power filtration captures the shape of the graph, just like Vietoris--Rips complexes capture the shape of a point cloud.
Power filtrations can be a crucial tool to employ in graph classification problems, and related questions.
As mentioned above, the higher persistence diagrams can contain crucial geometric information about graph cavities.
The main drawback here is that computing the power filtration for graphs and their higher-dimensional persistence diagrams is computationally expensive.
However, by using the ideas developed in this paper, there is hope to compute or bound these higher persistence diagrams more efficiently.

\bibliographystyle{plain}
\bibliography{ref.bib}

\appendix

\section{Proof of \cref{thm-main1}}
\label{app:proof}

Our proof of \cref{thm-main1} relies on Vietoris--Rips complexes of geodesic spaces and persistence modules that are indexed over the real numbers instead of over a discrete set, which we introduce now.
Our treatment is brief, since everywhere else in the paper we instead consider filtrations and persistence modules indexed over a finite index set.
However, we point the reader to references containing more information, including~\cite{bauer2014induced,gasparovic2018complete,virk20201}.

The Vietoris--Rips complex of a metric space is a filtration indexed over a real-valued scale parameter.

\begin{definition}
\label{def:VR}
For $X$ a metric space and $r\ge 0$, the \emph{Vietoris--Rips simplicial complex} $\vr{X}{r}$ has $X$ as its vertex set, and contains a finite subset $\sigma\subseteq X$ as a simplex if $\diam(\sigma)\le r$.
\end{definition}

Note that for any $r\le r'$, we have $\vr{X}{r}\subseteq \vr{X}{r'}$.
After applying homology $H_k$ with coefficients in $\Z/2\Z$, we obtain a persistence module indexed over the real numbers $r\ge 0$, namely the collection of vector spaces $\{H_k(\vr{X}{r})\}_{r\ge 0}$ equipped with linear maps $H_k(\vr{X}{r})\to H_k(\vr{X}{r'})$ for all $r\le r'$.
If $X$ is a compact metric space, then one can define the $k$-dimensional persistence diagram associated to the Vietoris--Rips filtration~\cite{ChazalDeSilvaOudot2014}, which (by a slight abuse of notation) we denote by $\PD_k(X)$.

For $G$ a finite connected graph, we let $\vr{G}{r}$ denote the Vietoris--Rips simplicial complex of the vertex set $\V$ of $G$, equipped with the shortest path metric $\rho$.
We have an equality of simplicial complexes $\vr{G}{r}=\wh{G}^{\lfloor r\rfloor}$ for any $r\ge 0$, where $\lfloor r \rfloor$ is the greatest integer smaller than or equal to $r$.

Let $\oG$ be the metric graph associated to the finite connected graph $G$.
To obtain $\oG$, start with the vertex set $\V$ of $G$, and then glue on a continuous interval $[0,1]$ of length $1$ for each edge in $G$.
One can equip with $\oG$ with a metric structure such that the distance between any two points (not necessarily vertices) is the length of a shortest path between them.
In particular, $\oG$ is a geodesic metric space that contains an infinite number of points on each edge of $G$, and that admits an isometric embedding from the vertex set of $G$; see~\cite{burago2001course} and~\cite[Section~1.9]{bridson2011metric} for more information.
Since $G$ is finite, the metric space $\oG$ is compact, and therefore for any integer $k$ we have a $k$-dimensional persistent homology diagram $\PD_k(\oG)$ coming from the Vietoris--Rips filtration of $\oG$~\cite{ChazalDeSilvaOudot2014}.

\begin{proof}[Proof of \cref{thm-main1}]
Let $G$ be a finite connected graph.
Let $\{\gamma_1,\ldots,\gamma_m\}$ be a lexicographically shortest basis of $H_1(G)$, where each $\gamma_i$ has length $l_i$.
So $l_1\le l_2\le l_3\le \ldots\le l_m$.
We must show that
\[\PD_1(G)=\{(\, 1,\left\lceil \tfrac{l_i}{3}\right\rceil\, )\mid 1\le i\le m\}.\]
By \cref{lem:birth-1}, we know that $\PD_1(G)=\{(1, d_i)\}_i$ for some collection of death times $d_i$.
We must show that $d_i=\lceil \frac{l_i}{3}\rceil$ for $1\le i\le m$.

The high-level structure of the proof will be as follows.
The persistence diagram of the metric graph $\PD_1(\oG)$ is known, by~\cite[Theorem~8.10]{virk20201}.
(We refer the reader to~ \cite[Theorem~1.1]{gasparovic2018complete} for an analogous result with \v{C}ech complexes.)
The papers~\cite{virk20201,gasparovic2018complete} are the inspiration for \cref{thm-main1}, and we rely on their results in this proof.
We will use morphisms of persistence modules~\cite{bauer2014induced} to pass from knowledge of the known persistence diagram $\PD_1(\oG)$ to knowledge of the unknown diagram $\PD_1(G)$.

A loop $\gamma$ in the finite graph $G$ is a sequence of adjacent edges $v_0v_1$, $v_1v_2$, \ldots, $v_{n-2}v_{n-1}$, $v_{n-1}v_0$ in $G$ that start and end at the vertex $v_0$.
By replacing each discrete edge $v_i v_{i+1}$ in the finite graph $G$ with the continuous path across this edge in the metric graph $\oG$, we obtain an associated loop $\ogamma$ of the same length in the metric graph $\oG$.
Note that $G$ is a $1$-dimensional simplicial complex triangulation of the metric space $\oG$.
By the equivalence between simplicial homology and singular homology, the rank of $H_1(G)$ is equal to the rank of $H_1(\oG)$.

Furthermore, the equivalence between simplicial homology and singular homology gives equivalence of lexicographically shortest bases, as follows.
%For each $\ogamma_i$ a loop in $\oG$, we say that $\ogamma_1,\ldots,\ogamma_m$ is a {\em lexicographically shortest basis} of $H_1(\oG)$ if the non-decreasing sequence of lengths $l_1\le l_2\le l_3\le \ldots\le l_m$ is lexicographically smallest among all bases for $H_1(\oG)$.
If $\gamma_1,\ldots,\gamma_m$ is a lexicographically shortest basis for $H_1(G)$, then $\ogamma_1,\ldots,\ogamma_m$ is a lexicographically shortest basis for $H_1(\oG)$ with the same lengths $l_1\le l_2\le l_3\le \ldots\le l_m$.
Theorem~8.10 of~\cite{virk20201} gives a complete description of the 1-dimensional persistent homology: the 1-dimensional persistence diagram of $\vr{\oG}{r}$ is given by $\PD_1(\oG)=\{(0,\frac{l_i}{3})\mid 1\le i\le m\}$, where the lengths $l_i$ are from a lexicographically shortest basis for $\oG$.

For any $r\le r'$, we have a commutative diagram
\begin{center}
\begin{tikzcd}
\vr{G}{r} \arrow[r, hook] \arrow[d, hook] & \vr{G}{r'} \arrow[d, hook] \\
\vr{\oG}{r} \arrow[r, hook] & \vr{\oG}{r'}
\end{tikzcd}
\end{center}
where the vertical maps are induced by the isometric inclusion of metric spaces from the vertex set of $G$ into the metric space $\oG$.
This means that we have a morphism of persistence modules, $f\colon H_1(\vr{G}{-})\to H_1(\vr{\oG}{-})$; see~\cite{bauer2014induced} for more background on morphisms between persistence modules.
The main result we will need is \cite[Proposition~5.3]{bauer2014induced}, which says that if $f$ is a morphism of persistence modules and if $(b,d)$ is a point in the persistence diagram of the domain of $f$ whose generator (at any scale between $b$ and $d$) gets mapped to the generator for a point $(b',d')$ in the persistence diagram of the codomain of $f$, then we have the inequality $b'\le b<d'\le d$.
In our context, this means the following, with $(b',d')=(0,\frac{l_i}{3})$ and with $b=1$.
Since the loop $\gamma_i$ in $G$ maps under the inclusion $\vr{G}{r}\hookrightarrow \vr{\oG}{r}$ to the homology class generated by the loop $\ogamma_i$ (for any $1\le r< l_i$), and since $\ogamma_i$ corresponds to the point $(0,\frac{l_i}{3})\in\PD_1(\oG)$, then $\gamma_i$ generates a bar $(1,d_i)\in\PD_1(G)$ that satisfies $0\le 1<\frac{l_i}{3}\le d_i$, for $1\le i\le m$.
As each $d_i$ is an integer, we furthermore have $d_i \ge \lceil\frac{l_i}{3}\rceil$.
By \cref{lem:death-3}, the cycle $\gamma_i$ is null-homotopic in $\wh{G}^n$ for $n\ge \lceil\frac{l}{3}\rceil$, which implies $d_i\le \lceil\frac{l_i}{3}\rceil$.
Together, these give $d_i=\lceil\frac{l_i}{3}\rceil$, and so $\PD_1(G)=\{(1, \lceil\frac{l_i}{3}\rceil)\mid 1\le i\le m\}$.
\end{proof}

\end{document}